\documentstyle[amssymb,11pt]{article}

\begin{document}

\baselineskip 15pt
\parindent=1em
\hsize=12.3 cm \textwidth=12.3 cm
\vsize=18.5 cm \textwidth=18.5 cm

\title{Understanding preservation theorems: Chapter VI of {\em Proper and Improper Forcing,} I}
\author{
Chaz Schlindwein \\
Department of Mathematics and Computing \\
Lander University \\
Greenwood, South Carolina 29649 USA\\
{\tt cschlind@lander.edu}}

\maketitle

\def\jux{\mathbin{\widehat{\ }}}
\def\forces{\mathbin{\Vdash}}

\def\restr{\mathbin{\upharpoonright}}

   \centerline{{\bf Abstract}}
   We present an exposition of Section VI.1 and most of Section VI.2 from Shelah's book {\em Proper and Improper Forcing}.  These sections offer proofs of the preservation under countable support iteration of proper forcing of various properties, including proofs that $\omega^\omega$-bounding, the Sacks property, the Laver property, and the $P$-point property are preserved by countable support iteration of proper forcing.  Also, any countable support iteration of proper forcing
that does not add a dominating real preserves ``no Cohen reals.''

\eject

\section{Introduction}

This paper is an exposition of some preservation theorems, due to Shelah [13, Chapter VI], for countable support iterations of proper forcing. 
These include the preservation of the ${}^\omega\omega$-bounding property, the Sacks and Laver properties, the $P$-point property, and some others.  Generalizations to revised countable support iterations of semi-proper forcings or even certain non-semi-proper forcings are given in [13, Chapter VI] but we do not address these more general iterations.  The results of [13, Section VI.2] overlap the results of [3] and [4], but the methods are dissimilar.  The article [1] covers similar ground.

This is the third in a sequence of expository papers covering parts of Shelah's book, {\it Proper and Improper Forcing.}  The earlier papers were [11], which covers sections 2 through 8 of [13, Chapter XI] and [9], which covers sections 2 and 3 of [13, Chapter XV]. 
The fourth paper of this sequence is [12], which presents  an exposition of
[13, Sections VI.3 and XVIII.3], including a proof of [13, Conclusion VI.2.15D].
 Other papers by the author generalize certain other results in [13]; in no instance were we content to quote a result of Shelah without supplying a proof.  Thus, [6] may be read, in part, as an exposition of [13, Sections V.6, IX.2, and IX.4]; [7] is, in part, an exposition of [13, Section V.8 and Theorem III.8.5];
and [8] includes as a special case an alternative proof of [13, Theorem III.8.6].  Also, [6] answers
[13, Question IX.4.9(1)];  [7] answers a question implicit in [13, Section IX.4]; [10] answers another such question and also may be read, in part, as an exposition of the results of Eisworth and Shelah [2] that weaken the assumption
$``\alpha$-proper for every $\alpha<\omega_1$'' used in [13, Section V.6]. 

\subsection{Notation}

We write $p\leq q$ when $p$ is a stronger forcing condition than $q$.

When $\langle P_\eta\,\colon\eta\leq\beta\rangle$ is a forcing iterartion, and $\alpha<\beta$, we set $P_\beta/{G_{P_\alpha}}$ to be a $P_\alpha$-name characterized by

\medskip

\centerline{$V[G_{P_\alpha}]\models``P_\beta/{G_{P_\alpha}}=\{p\restr[\alpha,\beta)\,\colon p\in P_\beta$ and $ p\restr\alpha\in G_{P_\alpha}\}.$''}

\medskip

Notice that when $p\in P_\alpha$ and $p\forces``q\in P_\beta/{G_{P_\alpha}}$'' then $p\forces``q\forces`\varphi$'\thinspace'' makes sense.  In contrast,
$p\jux q\forces``\varphi$'' makes sense only under the stronger assumption that $p\jux q\in P_\beta$.  For example, it could be the case that ${\rm supp}(p\jux q)=\beta$ yet
$p\forces``{\rm supp}(q)$ is a singleton.''
For this reason we favor the former notation and eschew the latter.

\section{Preservation of properness}

The fact that properness is preserved under countable support iterations was proved by Shelah in 1978.  The proof of this fact 
 is the basis of all preservation theorems for countable support iterations.

\proclaim Theorem 2.1 (Proper Iteration Lemma, Shelah). Suppose $\langle P_\eta\,
\colon\allowbreak\eta\leq\kappa\rangle$
is a countable support forcing iteration based on\/
$\langle Q_\eta\,\colon\allowbreak\eta<\kappa\rangle$ and
for every $\eta<\kappa$ we have that\/ ${\bf 1}\forces_{P_\eta}`` Q_\eta$ {\rm is
proper.''} Suppose also that\/ $\alpha<\kappa$ and\/
$\lambda$ is a sufficiently large regular cardinal and\/ $N$ is a
countable elementary submodel of\/ $H_\lambda$ and\/
$\{P_\kappa,\alpha\}\in N$ and\/
$p\in P_\alpha$ is\/ $N$-generic and $q$ is a $P_\alpha$-name  and\/ {\rm
$p\forces``q\in  P_\kappa/G_{P_\alpha}\cap
N[G_{P_\alpha}]$.''}
Then there is $r\in P_\kappa$ such that
$r$ is $N$-generic and $r\restr\alpha=p$ and\/ {\rm
$p\forces``r\restr[\alpha,\kappa)\leq
q$.''}

Proof.  The proof proceeds by induction, so suppose that
the Theorem holds for all iterations of length less than $\kappa$.  Fix $\lambda$  and  $N$  and  $\alpha$ and $p$ and
 $q$ as in the assumption.

Case 1.  $\kappa=\beta+1$ for some $\beta$.

Because $\beta\in N$  we may use the induction hypothesis to fix $p'\in P_\beta$ such that
$p'\restr\alpha=p$ and $p'$ is $N$-generic and
$p\forces``p'\restr[\alpha,\beta)\leq q\restr\beta$.''  We have that
$p'\forces``q(\beta)\in N[G_{P_\beta}]$.'' Take $r\in P_\kappa$ such that
$r\restr\beta=p'$ and 

\medskip

{\centerline{$p'\forces``r(\beta)\leq q(\beta)$ and $r(\beta)$ is
$N[G_{P_\beta}]$-generic for $Q_\beta$.''}}

\medskip

  Then $r$ is $N$-generic and we are done with the successor case.

Case 2.  $\kappa$ is a limit ordinal.

Let $\beta={\rm sup}(\kappa\cap N)$, and fix $\langle\alpha_n\,\colon\allowbreak
n\in\omega\rangle$ an increasing sequence from $\kappa\cap N$ cofinal in $\beta$ such that
$\alpha_0=\alpha$.  Let $\langle\sigma_n\,\colon\allowbreak n\in\omega\rangle$ enumerate all the $P_\kappa$ names $\sigma\in N$ such that ${\bf 1}\forces``\sigma$ is an ordinal.''

Using the induction hypothesis, build a sequence $\langle \langle p_n,q_n,\tau_n\rangle\,\colon\allowbreak n\in\omega\rangle$ such that $p_0=p$ and $p\forces``q_0\leq q$'' and for each $n\in\omega$ we have all of the following:

(1) $p_n\in P_{\alpha_n}$ and $p_n$ is $N$-generic  and
$p_{n+1}\restr\alpha_n=p_n$ and $\tau_n$ is a $P_{\alpha_n}$-name.

(2)  $p_n\forces``q_n\in P_{\kappa}/G_{P_{\alpha_n}}\cap N[G_{P_{\alpha_n}}]$ and
 $\tau_{n}\in  N[G_{P_{\alpha_n}}]  $ and $q_n\forces`\sigma_{n}=\tau_{n}
$' and
if $n>0$ then $q_n\leq q_{n-1}\restr[\alpha_n,\kappa)$.''

(3) $p_n\forces``p_{n+1}\restr[\alpha_n,\alpha_{n+1})\leq q_n\restr\alpha_{n+1}$.''

Define $r\in P_\kappa$ such that $(\forall n\in\omega)\allowbreak
(r\restr\alpha_n=p_n)$ and ${\rm supp}(r)\subseteq\beta$.  
To see that $r$ is $N$-generic, suppose that $\sigma\in N$ is a $P_\kappa$-name for an ordinal. Fix $n$ such that
$\sigma = \sigma_n$. 
Because $p_{n}$ is $N$-generic, we have 

\medskip

\centerline{$p_{n}\forces``{\rm supp}(q_{n})\subseteq\kappa\cap N[G_{P_{\alpha_{n}}}]=
\kappa\cap N$,''}

\medskip

\noindent whence it is clear that 

\medskip

\centerline{$p_{n}\forces``r\restr[\alpha_{n},\kappa)\leq 
q_{n}$.'' }

\medskip

We have 

\medskip

\centerline{$p_{n}\forces``q_{n}\forces`\sigma\in Ord\cap N[G_{P_{\alpha_{n}}}]=Ord 
\cap N
$,'\thinspace''}

\medskip

\noindent where $Ord$ is the class of all ordinals. 
Thus $r\forces``\sigma\in N$.''
We conclude that $r$ is $N$-generic, and the Theorem is established.

\proclaim Corollary 2.2 (Fundamental Theorem of Proper Forcing, Shelah).
Suppose $\langle P_\eta\,
\colon\allowbreak\eta\leq\kappa\rangle$
is a countable support forcing iteration based on\/
$\langle Q_\eta\,\colon\allowbreak\eta<\kappa\rangle$ and
for every $\eta<\kappa$ we have that\/ ${\bf 1}\forces_{P_\eta}`` Q_\eta$ {\rm is
proper.''}  Then $P_\kappa$ is proper.

Proof: Take $\alpha = 0$ in the Proper Iteration Lemma.

\section{Preservation of proper plus $\omega^\omega$-bounding}

In this section we recount Shelah's proof of the preservation of ``proper plus $\omega^\omega$-bounding.''  This is a special case of [13, Theorem VI.1.12] and is given as [13, Conclusion VI.2.8D].
Another treatment of this result can be found in [3] and [4], using different methods.
Notably, the proof given in [3] assumes each forcing adds reals and [4] contains a patch for this deficiency, but in
Shelah's proof presented here, the issue does not arise.

The following Lemma justifies the construction of $\langle p^n\,\colon n\in\omega\rangle$
and $\langle t_{n,m}\,\colon\allowbreak m\leq n<\omega\rangle$ 
and $\langle f_m\,\colon\allowbreak
m\in\omega\rangle$ in [13, proof of Theorem VI.1.12], where Shelah's $p^n$ is our $p_n\restr n$.
Our $p'$ encapsulates the third paragraph of
[13, proof of Theorem VI.i.12], i.e., Shelah's assertion that w.l.o.g.~$f(k)$ is a $P_k$-name of
a natural number (see (2) below).

\proclaim Lemma 3.1.  Suppose $\langle P_n\,\colon\allowbreak
n\leq\omega\rangle$ is a countable support iteration 
based on $\langle Q_n\,\colon\allowbreak n<\omega\rangle$.
Suppose also that $f$\/ is a $P_\omega$-name for an element
of\/ ${}^\omega\omega$, and suppose $p\in P_\omega$.
Then there are
$\langle p_n\,\colon\allowbreak n\in\omega\rangle$ and $\langle f_n\,\colon\allowbreak n\in\omega\rangle$
and $p'\leq p$  such that
$p_0= p'$ and for every $n\in\omega$ we have that each of the following holds:

(1) $f_n$ is a $P_{n}$-name for an element of ${}^\omega\omega$, and

(2) $p'\restr n\forces``p'\restr[n,\omega)\forces`f\restr(n+1)=f_n\restr(n+1)$,'\thinspace'' and

(3) $p_{n+1}\leq p_n$, and

(4) for all $m$, we have $p'\restr n\forces``p_m(n)
\forces`f_n\restr(m+1)=f_{n+1}\restr(m+1)$.'\thinspace''

\medskip

Proof: Define $\langle q_n\,\colon n\in\omega\rangle$ and $\langle\sigma_n\,\colon\allowbreak n\in\omega\rangle$ such that $\sigma_0\in\omega$ and $q_0\leq p$ and  $q_0\forces_{P_\omega}``f(0)=\sigma_0$''
and for every $n>0$ we have that $\sigma_n$ is a $P_n$-name for an integer and
${p\restr n}\forces_{P_{n}}``q_n\in P_{\omega}/G_{P_n}$ and $q_n\leq q_{n-1}\restr[n,\omega)$ and
$q_n\forces`\sigma_n=f(n)$.'\thinspace''

Define $p'\in P_\omega$ by $(\forall n\in\omega)(p'(n)=q_n(n))$.
We have that (2) holds.

We now define $\langle p_n\,\colon n>0\rangle$.  Given $p_n$, build $\langle( q_n^k
,\tau^k_n)\,\colon\allowbreak k\leq n\rangle$ by downward induction by setting $\tau_n^{n+1} =\sigma_{n+1}$ and, given $\tau_n^{k+1}$, take 
$q_n^k$ and $\tau_n^k$ such that $\tau_n^k$ is a $P_k$-name for an integer
and $p'\restr 
k\forces``q_n^k\leq
 p_{n}(k)$ and 
$q_n^k\forces`\tau_n^k=\tau_n^{k+1}$.'\thinspace''
Choose $p_{n+1}\leq p_n$ such that for every $k\leq n$ we have $p_{n+1}(k)= q_n^k$.

For every $k\leq n$ let $f_n(k)$ be the $P_{n}$-name for $\sigma_k$ and for $k>n$ let $f_n(k)=\tau_k^n$.

The Lemma is established.

\proclaim Definition 3.2.  For $f$ and $g$ in ${}^\omega\omega$ we say
$f\leq g$ iff $(\forall n\in\omega)\allowbreak(f(n)\leq g(n))$.
We say that $P$ is ${}^\omega\omega$-bounding iff\/ 
$V[G_P]\models``(\forall f\in{}^\omega\omega)\allowbreak(\exists g\in {}^\omega\omega\cap V)
\allowbreak(f\leq g)$.''

\proclaim Lemma 3.3. Suppose $P$ is ${}^\omega\omega$-bounding and\/ {\rm $V[G_P]\models``(\forall n\in\omega)
\allowbreak(f_n\in {}^\omega\omega)$.''}
Then\/ {\rm $V[G_P]\models``(\exists\langle f'_n\,\colon\allowbreak n\in\omega\rangle\in V)\allowbreak
(\forall n\in\omega)\allowbreak(f'_n\in {}^\omega\omega$ and $f_n\leq f'_n)$.''}

Proof: Let $j$ be a one-to-one function from ${}^2\omega$ onto $\omega$. In $V[G_P]$ define $h\in {}^\omega\omega$ by $(\forall n\in\omega)\allowbreak(\forall m\in\omega)\allowbreak(h(j(n,m))=f_n(m))$ and
choose $h'\in {}^\omega\omega\cap V$ such that $h\leq h'$.  Define $\langle f'_n\,\colon\allowbreak
n\in\omega\rangle$ by $(\forall n\in\omega)\allowbreak(\forall m\in\omega)\allowbreak(h'(j(n,m))=f'_n(m))$.
The Lemma is established.

The proof  of the following Theorem is obtained from [13, proof of Theorem VI.1.12] by discarding all references to $x$ and replacing each tree with a function bounding its branches; both of these simplifications are justified by [12, Definition VI.2.8A].  Also we have removed any reference to non-proper forcings and we have made explicit the dependence of the
 functions $F_0$, $F_1$, and $F_2$ on the parameter $n$ (this dependence is suppressed in Shelah's proof).

\proclaim Theorem 3.4.  Suppose $\langle P_\eta\,\colon\allowbreak\eta\leq\kappa\rangle$ is a countable support iteration based on $\langle Q_\eta\,\colon\allowbreak\eta<\kappa\rangle$ and suppose\/ {\rm $(\forall\eta<\kappa)\allowbreak({\bf 1}\forces_{P_\eta}``Q_\eta$ 
is proper and ${}^\omega\omega$-bounding'').} Then $P_\kappa$ is ${}^\omega\omega$-bounding.

Proof:  By induction on $\kappa$.  By 
standard arguments, taking into account the fact that a
 counterexample to ${}^\omega\omega$-bounding
cannot first appear in $V[G_{P_\kappa}]$ where $\kappa$ has uncountable cofinality,
and the fact that the composition of two ${}^\omega\omega$-bounding 
forcings is again ${}^\omega\omega$-bounding,
 we only have to establish this for $\kappa=\omega$.

Fix $\lambda$ a sufficiently large regular cardinal and $N$ a countable elementary substructure
of $H_\lambda$ such that $P_\omega\in N$ and suppose $ p\in P_\omega\cap N$.

Let $\langle g_j\,\colon\allowbreak j<\omega\rangle$ list 
${}^\omega\omega\cap N$, with infinitely many repetitions.

Fix $p'$ and $\langle (p_n,f_n)\,\colon\allowbreak n\in\omega\rangle$ as in
Lemma 3.1.
We may assume that for every $n\in\omega$ we have $p'$ and $p_n$ and $f_n$
are in $N$, and, furthermore, the sequence  $\langle\langle 
p_n,f_n\rangle\,\colon\allowbreak n\in\omega\rangle$ is in $ N$.

Define $ g\in{}^\omega\omega$ by
\medskip

{\centerline{$g(i)={\rm max}\{f_0(i),\allowbreak{\rm max}\{g_k(i)\,\colon\allowbreak
k\leq i\}\}$.}}

\medskip

For each $n\in\omega$, fix $P_{n}$-names $F_{n,0}$ and $F_{n,2}$
such that $V[G_{P_{n}}]\models``F_{n,0}$ maps $Q_n$ 
into ${}^\omega\omega$ and
$F_{n,2}$ maps $Q_n$ into $Q_n$ and
for every $q'\in Q_n$ we have
$F_{n,2}(q')\leq q'$ and $F_{n,2}(q')\forces`f_{n+1}\leq F_{n,0}(q')$'.''

For each $n\in\omega$, use Lemma 3.3 to fix $F_{n,1}$ such that $V[G_{P_{n}}]\models``F_{n,1}\in
V$ maps $\omega$ to ${}^\omega\omega$ and for all $m\in\omega$ we have
$F_{n,0}(p_m(n))\leq F_{n,1}(m)$.''

We may assume that for each $n\in\omega$ we have the name $F_{n,1}$ is in $N$.

Claim 1.  We may be build $\langle r_n\,\colon\allowbreak n\in\omega\rangle$ such that
 for every $n\in\omega$ we have that the following hold:

(1) $r_n\in P_{n}$ is $N$-generic.

(2) $r_{n+1}\restr n=r_n$.

(3) $r_n\forces``f_n\leq g$.''

(4) $r_n\forces``r_{n+1}(n)\leq p'(n)$.''

Proof: Work by induction on $n$.  

Case 1:  $n=0$.

We have $f_0\leq  g$.

Case 2: Otherwise.

Suppose we have $r_n$.

In $V[G_{P_{n}}]$, define $ g_n^*$ by
$(\forall i\in\omega)(g_n^*(i)= 
{\rm max}\{F_{n,1}(m)(i)\,\colon\allowbreak m\leq i\})$.

We may assume the name $g_n^*$ is in $N$.

Notice that we have

\medskip

\centerline{$r_n\forces``g_n^*\in N[G_{P_{n}}]\cap V = N
$.''}

\medskip

Therefore we may choose a $P_{n}$-name $k$ such that $r_n\forces``g_n^*=g_k$ and $k>n$'' (in our notation 
we suppress the fact that $k$ depends on $n$).

Subclaim 1: $r_n\forces``F_{n,2}(p_k(n))
\forces`f_{n+1}\leq g$.'\thinspace''

Proof: 
For $i\geq k$ we have 

\medskip

{\centerline{$r_n\forces``F_{n,2}(p_k(n))
\forces`f_{n+1}(i)\leq F_{n,0}(p_k(n))(i)$}}

{\centerline{${}\leq F_{n,1}(k)(i)\leq g_n^*(i)=g_k(i)\leq g(i)$.'\thinspace''}}

\medskip

\noindent The first inequality is by the definition of $F_{n,0}$, the second inequality is
by the definition of $F_{n,1}$,
the third inequality is by the definition of $g_n^*$ 
along with the fact that $i\geq k$, the  equality is by the 
definition of $k$, and the last inequality is by the definition 
of $g$ along with the fact that $i\geq k$.

For $i<k$, we have 

\medskip

\centerline{$r_n\forces``p_k(n)
\forces`f_{n+1}(i)=f_n(i)\leq g(i)$.'\thinspace''}

\medskip

\noindent The equality is by the choice of $\langle (f_m, p_m)\,\colon m\in\omega\rangle$ (see Lemma 3.1),
and the inequality is by the induction hypothesis that Claim 1  holds for integers less than or equal to $n$.

Because $r_n\forces``F_{n,2}(p_k(n))\leq p_k(n)$,'' we have that the Subclaim is established.

Choose $r_{n+1}\in P_{{n+1}}$ such that
$r_{n+1}$ is $N$-generic and $r_{n+1}\restr n=r_n$ and

\medskip

\centerline{
$r_n\forces``r_{n+1}(n)\leq F_{n,2}(p_k(n))$.''}

\medskip

This completes the proof of Claim 1.

Let $r=\bigcup\{r_n\,\colon\allowbreak n\in\omega\}$.  We have that

\medskip

\centerline{
$r\leq p$ and 
$r\forces``f\leq g$.''}

\medskip

The Theorem is established.

\section{The Sacks property}

In this section we present Shelah's proof of the preservation of ``proper plus Sacks property'' under countable support iteration.  The proof is a special case of  [13, Theorem VI.1.12]
 and appears as [13, Conclusion VI.2.9D].  

\proclaim Definition 4.1.  For $x$ and $y$ in ${}^\omega(\omega-\{0\})$,
 we say that $x\ll y$ iff $(\forall n\in\omega)\allowbreak(x(n)\leq y(n))$ and
\[ \lim_{n\rightarrow\infty}y(n)/x(n)=\infty \]
In particular for $x\in{}^\omega(\omega-\{0\})$ we have $1\ll x$ iff\/ $x$ diverges to infinity.

The following Definition corresponds to [13, Definition VI.2.9A(b)].

\proclaim Definition 4.2.  
For\/ $T\subseteq{}^{<\omega}\omega$ a tree and\/ $x\in{}^\omega(\omega-\{0\})$, 
we say that $T$ is an $x$-sized tree iff for every $n\in\omega$ we have that the 
cardinality of\/ $T\cap{}^n\omega$ is at most $x(n)$, and\/ $T$ has  no terminal nodes.

\proclaim Definition 4.3.  For\/ $T\subseteq{}^{<\omega}\omega$ we set\/ $[T]$
equal to the set of all\/ $f\in{}^\omega\omega$ such that
every initial segment of\/ $f$ is in $T$.  That is, $[T]$ is the set of infinite branches of\/ $T$.

Often the Sacks property is given as a property
of pairs of models; however, because our focus is on forcing constructions, we define it to be a property of posets.

\proclaim Definition 4.4.  
A poset\/ $P$ has the Sacks property iff whenever\/ $x\in{}^\omega(\omega-\{0\})$ and\/ 
$1\ll x$
 then we have

{\centerline{${\bf 1}\forces_P``(\forall f\in{}^\omega\omega)\allowbreak(\exists H\in V)\allowbreak(H$ is an $x$-sized tree and $f\in[ H])$.''}}

\proclaim Definition 4.5.  Suppose $n\in\omega$. We say that $t$ is an $n$-tree iff $t\subseteq{}^{\leq n}\omega$ and $t$ is closed under initial segments and $t$ is non-empty and for every $\eta\in t$ there is $\nu\in t$ such that $\nu $ extends $\eta$ and ${\rm lh}(\nu)=n$.

The following Lemma shows that $(D,R)$ given in [13, Definition VI.2.9A] satisfies
[13, Definition VI.2.2(3)$(\varepsilon)^+$].  The proof given here follows [13, proof of Claim VI.2.9B$(\varepsilon)^+$].

\proclaim Lemma 4.6.  $P$ has the Sacks property iff whenever  $x$ and  $z$ are elements of\/
 ${}^\omega(\omega-\{0\})$ and $x\ll z$ we have that

{\centerline{${\bf 1}\forces_P``(\forall T)($if $T$ is an $x$-sized tree then}}
{\centerline{$(\exists H\in V)\allowbreak(H$ is a $z$-sized tree and $T\subseteq H))$.''}}

\medskip

Proof: We prove the non-trivial direction.
Suppose $P$ has the Sacks property and suppose $x$ and $z$ are given.   Working in $V[G_P]$, 
suppose $T$ is given.  For every $n\in\omega$ let 

\medskip

\centerline{${\cal T}_n(x)=\{t\subseteq{}^{\leq n}\omega\,\colon\allowbreak
t$ is an $n$-tree}

\centerline{and $(\forall i\leq n)\allowbreak(\vert t\cap {}^i\omega\vert\leq x(i))\}$.}

\medskip

Let

\medskip

\centerline{${\cal T}(x)=\bigcup\{{\cal T}_n(x)\,\colon\allowbreak n\in\omega\}$.}

\medskip

 Under the natural order, 
${\cal T}(x)$ is isomorphic to ${}^{<\omega}\omega$.

Define $\zeta\in[{\cal T}(x)]$ by setting $\zeta(n)=T\cap {}^{\leq n}\omega$ for all $n\in\omega$.

Define $y\in {}^\omega(\omega-\{0\})$ by setting $y(n)$ equal to the greatest integer less than or equal to $z(n)/x(n)$ for every $n\in\omega$.
Clearly $1\ll y$, so we may choose a $y$-sized tree $H'\subseteq {\cal T}(x)$ such that $\zeta\in[H']$ and
$H'\in V$.

Let $H=\bigcup H'$.
We have that $H$ is a $z$-sized tree and $H\in V$ and $T\subseteq H$.

The Lemma is established.

The  following Lemma is 
[13, Claim VI.2.4(1)] specialized to the case of 
Sacks property, and we follow the proof from [13].

\proclaim Lemma 4.7. Suppose $y$ and $z$ are elements of ${}^\omega(\omega-\{0\})$
and $y\ll z$. Suppose $P$ is a forcing such that\/ {\rm
$V[G_P]\models``$for every countable $X\subseteq V$ there is a countable $Y\in V$ such that $X\subseteq Y$.''}
Suppose in $V[G_P]$ we have that\/  $\langle T_n\,\colon\allowbreak n\in\omega\rangle$
is a sequence of $y$-sized trees such that for every $n$ we have $T_n\in V $.
Then in $V[G_P]$ there is a $z$-sized tree $T^*\in V$ and an increasing sequence of integers $\langle m(n)\,\colon  n\in\omega\rangle$
such that\/ $m(0)=0$ and\/ 
$(\forall n> 0)\allowbreak(m(n)>n)$ and for every $\eta\in {}^{<\omega}\omega
$ we have

\medskip

\centerline{$(\forall n> 0)\allowbreak(\exists i<n)(\eta\restr
 m({n+1})\in T_{m(i)})$ implies $\eta\in T^*$.}

\medskip

Proof: Fix $x\in {}^\omega(\omega-\{0\})$ such that $y\ll x\ll z$. Fix $\langle x_n\,\colon\allowbreak
n\in\omega\rangle$ a sequence of elements of ${}^\omega(\omega-\{0\})$ such that
$(\forall n\in\omega)(y\ll x_n\ll x_{n+1}\ll x)$.

Set $k(0)=0$ and for each $n>0$ set $k(n)$ equal to the least $k>n$ such that

\medskip

\centerline{$(\forall j\geq k)(2x_n(j)\leq x_{n+1}(j)$ and $(n+1)x(j)\leq z(j))$.}

\medskip

Work in $V[G_P]$. Let $b\in V$ be a countable set of $y$-sized trees such that $\{T_n\,\colon\allowbreak n\in\omega\}\subseteq b$.
Let $\langle S_n\,\colon\allowbreak n\in\omega\rangle\in V$ enumerate $b$ with infinitely 
many repetitions with $S_0=T_0$.

Build $\langle S'_n\,\colon n\in\omega\rangle$ by setting $S'_0=S_0$ and for every $n\in\omega$ let

\medskip

(A)  $S'_{n+1}=\{\rho\in S_n\,\colon\allowbreak\rho\restr k(n)\in S'_n\}\cup S'_n$.

\medskip

Claim 1. For all $n\in\omega$ we have that $S'_n$ is an $x_n$-sized tree.

Proof: By induction on $n$. Clearly $S'_0=T_0$ is an $x_0$-sized tree. For every $t<k(n)$ we have that
$\vert S'_{n+1}\cap{}^t\omega\vert= \vert S'_n\cap{}^t\omega\vert\leq x_n(t)\leq x_{n+1}(t)$.
For every $t\geq k(n)$ we have $\vert S'_{n+1}\cap{}^t\omega\vert\leq
\vert S'_n\cap{}^t\omega\vert+\vert S_n\cap{}^t\omega\vert\leq x_n(t)+y(t)\leq x_{n+1}(t)$. 
The Claim is established.

Define $h\in{}^\omega\omega$ by setting $h(0)=0$ and for every $n>0$ setting
$h(n)$ equal to the least $m>n$ such that $T_n=S_m$.

Build $\langle n'_i\,\colon\allowbreak i\in\omega\rangle$ an increasing sequence of integers such that
$n'_0=0$ and $n'_1>k(1)$ and for every $i\in\omega$ we have

(B) $k(h(n'_i))<n'_{i+1}$.

By (B) we have

(C) $(\forall i\in\omega)(\exists t)(n'_i<k(t)<n'_{i+1})$.

Let $T^*=\{\eta\in{}^{<\omega}\omega\,\colon\allowbreak (\forall n>0)(\exists i<n)\allowbreak(\eta\restr k(n)\in S'_{k(i)}
)\}$.

Claim 2. $T^*$ is a $z$-sized tree.

Proof.  Given $t\geq k(1)$,  choose $n\in\omega$ such that $k(n)\leq t<k({n+1})$.
We have 

\medskip

\centerline{$ T^*\cap {}^t\omega\subseteq \{\eta\in{}^t\omega\,\colon\allowbreak(\forall j\leq n+1)\allowbreak
($if $j>0$ then $(\exists i< j)\allowbreak(\eta\restr k(j)\in S'_{k(i)}))\}$}

\medskip

\noindent and so 

\medskip

\centerline{$\vert T^*\cap {}^t\omega\vert\leq\Sigma_{i\leq n}\vert S'_{k(i)}\cap {}^t\omega
\vert\leq(n+1)x(t)
\leq z(t)$.}

\medskip

For $t<k(1)$ we have $T^*\cap{}^t\omega=T_0\cap{}^t\omega$, so $\vert T^*\cap{}^t\omega\vert\leq y(t)\leq z(t)$.

The Claim is established

For every $i\in\omega$ let $m_i=n'_{4i}$.

Fix $\eta\in{}^{<\omega}\omega$ such that

(D) $(\forall i> 0)(\exists j<i)(\eta\restr m_{i+1}\in T_{m_j})$.

To establish the Lemma, it suffices to show 

(E) $(\forall i>0)(\exists j<i)(\eta\restr k(i)\in S'_{k(j)})$,

\noindent since this implies $\eta\in T^*$.

Claim 3. $(\forall i>0)(\exists j<i)(\eta\restr n'_{i+1}\in S'_{n'_j})$.

We prove this by induction on $i$.

Case 1. $i<8$.

We have $n_{i+1}'\leq n'_8=m_2$ and by (D) we have $\eta\restr m_2\in T_0$.
Therefore $\eta\restr n'_{i+1}\in T_0=S_0=S'_0$.

Case 2. $i\geq 8$,

Fix $i^*$ such that $4i^*\leq i<4i^*+4$.

By (D) we may fix $j^*<i^*$ such that

(F) $\eta\restr m_{i^*+1}\in T_{m_{j^*}}=S_{h(m_{j^*})}$.

Using the fact that $4j^*+1<i$ we have, by the induction hypothesis, that

(G) $\eta\restr n'_{4j^*+1}\in S'_{n'_{4j^*}}=
S'_{m_{j^*}}\subseteq S'_{h(m_{j^*})}$.

By  (B) we have

(H) $k(h(m_{j^*}))<n'_{4j^*+1}$.

By (A), (F), (G), and (H) we have

(I) $\eta\restr m_{i^*+1}\in S'_{h(m_{j^*})+1}$.

We have

(J) $n'_{i+1}\leq n'_{4i^*+4}=m_{i^*+1}$.

By (I) and (J) we have $\eta\restr n'_{i+1}\in S'_{h(m_{j^*})+1}\subseteq S'_{n'_{4j^*+2}}$.

Because $4j^*+2<i$ the Claim is established.

To complete the proof of the Lemma, suppose $i>0$. By (E) it suffices to
show that there is $t<i$ such that
$\eta\restr k(i)\in S'_{k(t)}$.

Case 1: $k({i-1})<n'_0$.

By (C) we have $n'_1\geq k({i})$. By Claim 3 we have $\eta\restr n'_2\in S_0$.  
Hence $\eta\restr n'_1\in S_0$. Hence
$\eta\restr k(i)\in S_0$.

Case 2:  $n'_0\leq k({i-1})$.

By (C) we know that there is at most one element of $\{n'_j\,\colon\allowbreak j\in\omega\}$ strictly between
$k({i-1})$ and $k(i)$. Hence we may
fix $j>0$ such that $n'_{j-1}\leq k({i-1})<k({i})\leq n'_{j+1}$.  If $\eta\restr n'_{j+1}\in S_0$ then
$\eta\restr k({i})\in S_0$ and we are done, so assume otherwise.  By Claim 3 we may fix $m<j$ such that $\eta\restr n'_{j+1}\in S'_{n'_m}$.
We have $\eta\restr k(i)\in S'_{n'_m}\subseteq S'_{n'_{j-1}}\subseteq S'_{k({i-1})}$ and again we are done.

The Lemma is established.

The following Lemma shows that $(D,R)$ from [13, Definition VI.2.9A] satisfies [13, Definition VI.1.2$(\beta)(iv)$ and Remark VI.1.3(8)]
The proof given here is [13, proof of Claim VI.2.9B$(\gamma^+)$].

\proclaim Lemma 4.8. Suppose  $x\in{}^\omega(\omega-\{0\})$ 
and $z\in{}^\omega(\omega-\{0\})$ 
and  $x\ll z$.  
Suppose that $\langle T_n\,\colon\allowbreak n\in\omega\rangle$
 is a sequence of $x$-sized trees and $T$ is an $x$-sized tree.
Then  there is a $z$-sized tree $T^*$ and a sequence of integers
$\langle k_i\,\colon\allowbreak i\in\omega\rangle$
such that $T\subseteq T^*$ and
for every $\eta\in T$ and $i\in\omega$ and every $\nu\in T_{m_i}$ extending $\eta$, if 
${\rm length}(\eta)\geq k_i$  then  $\nu\in T^*$.

Proof. Choose $y\in{}^\omega(\omega-\{0\})$ such that $x\ll y\ll z$.

Fix $n^*$ such that $(\forall i\geq n^*)(2x(i)\leq y(i))$.

For every $n\geq n^*$ let $T'_n=\{\eta\in T_n\,\colon\allowbreak\eta\restr n\in T\}\cup T$, and for
$n<n^*$ let $T'_n=T$.

We have $T'_n$ is a $y$-sized tree for every $n\in\omega$.

For each $n\in\omega $ set $k_n$ equal to the least $k\geq n^*$ such that
$(\forall j\geq k)\allowbreak((n+2)y(j)\leq z(j))$.

Let $T^*=\{\eta\in{}^\omega\omega\,\colon\allowbreak(\forall i>0)\allowbreak(\exists j<i)\allowbreak
(\eta\restr k_i\in T'_{k_j})\}$.

Clearly $T\subseteq T^*$.

Claim: $T^* $ is a $z$-sized tree.

Proof: Like Claim 2 of the proof of Lemma 4.7.

Now suppose that $\eta\in T$ and $i\in\omega$ and
${\rm length}(\eta)\geq k_i$ and
 $\nu$ extends $\eta$
and $\nu\in T_{k_i}$.  We show $\nu\in T^*$.

Because $\nu$ extends an element of $T$ of length at least $k_i$, we have that $\nu\in T'_{k_i}$.
Choose $h\in[T'_{k_i}]$ extending $\nu$.  It suffices to show that $h\in[T^*]$.
Therefore it suffices to show that for every $n\in\omega$ we have

\begin{itemize}

\item[$(*)_n$]  $(\exists j< n)(h\restr k_n\in T'_{k_j})$.

\end{itemize}

Fix $n\in\omega$.

Case 1: $i<n$.

Because $h\in[T'_{k_{i}}]$ we have $h\restr k_n\in T'_{k_i}$, so $(*)_n$ holds.  

Case 2: $n\leq i$.

We have $h\restr k_n=\nu\restr k_n=\eta\restr k_n\in T\subseteq T'_{k_0}$.  Therefore $(*)_n$ holds.

The Lemma is established.

The proof of the following Lemma is ``proof of $(d)$'' in [13, proof of Claim VI.1.8].
Thus, it shows that $(D,R)$ from [13, Definition VI.2.9A] is a strong covering model [13,
Definition VI.1.6(4)].  Indeed, we (and Shelah)
show something stronger insofar as 
the quantifier ``there exists $\langle x_n\,\colon\allowbreak n\in\omega\rangle$'' 
in [13, Definition VI.1.6(4)(d)] is replaced  with the quantifier ``for all increasing $\langle x_n\,\colon\allowbreak n\in\omega\rangle$
bounded below $z$.''

\proclaim Lemma 4.9.  Suppose\/
$y\in{}^\omega(\omega-\{0\})$ and\/ $z\in{}^\omega(\omega-\{0\})$, and suppose
$\langle x_n\,\colon\allowbreak n\in\omega\rangle$ is a sequence of elements
of\/ ${}^\omega(\omega-\{0\})$ such that $(\forall n\in\omega)\allowbreak
(x_n\ll x_{n+1}\ll y\ll z)$.  
Suppose
 for every $n\in\omega$, we have $x_n^*\in{}^\omega(\omega-\{0\})$ and\/
$x_n^*\ll x_n\ll x^*_{n+1}$, and we have
$\langle x_{n,j}\,\colon\allowbreak n\in\omega$, $j\in\omega\rangle$ is a
sequence of elements of ${}^\omega(\omega-\{0\})$ such that
for every  $j\in\omega$ we have
$x_n\ll x_{n,j}\ll x_{n,j+1}\ll x_{n+1}^*$.
Suppose that $P$ is a forcing such that\/ {\rm $V[G_P]\models``$for every
countable $X\subseteq V$ there is a countable $Y\in V$ such that $X\subseteq Y$.''}
Suppose in $V[G_P]$ we have that $T\in V$ is an $x_0$-sized tree and
 $\langle T_{n,j}\,\colon\allowbreak n\in\omega$, $j\in\omega\rangle$ is
a sequence such that for every $n\in\omega$ and
 $j\in\omega$ we have that $T_{n,j}\in V$ is an
$x_{n,j}$-sized tree (but the sequence need not be in $V$).  Then in $V[G_P]$ there is a
sequence $\langle T^n\,\colon\allowbreak n\in\omega\rangle$ and
$T^*\in V$ such that $T^*$ is a $z$-sized tree and $T\subseteq T^*$ and for every
$n\in\omega$ we have

(i) $T^n\subseteq T^{n+1}$ and $T^n\in V$ is an $x_n$-sized tree, and

(ii)   for every $j\in\omega$ and every $g\in[T_{n,j}]$ 
there is $k\in\omega$ such that for every $\eta\in T_{n,j}$ extending $g\restr k$, if 
 $\eta\restr k\in T^n\cap T^*$ then
$\eta\in T^{n+1}\cap T^*$.

\medskip

Proof:  Work in $V[G_P]$.  Let $T^0=T$. 
Given $T^n\in V$, build $\langle T'_{n,j}\,\colon\allowbreak j\in\omega\rangle$
as follows.  Let $T'_{n,0}=T^n$.  Given
$T'_{n,j}$ take $m(n,j)\in\omega$ such that

\medskip

\centerline{
$(\forall t\geq m(n,j))\allowbreak(2x_{n,j}(t)\leq x_{n,j+1}(t))$.}

\medskip

Let $T'_{n,j+1}=\{\eta\in T_{n,j}\colon\allowbreak
\eta\restr m(n,j)\in T'_{n,j}\}\cup T'_{n,j}$.

Claim 1.  Whenever $i\leq j<\omega$ we have $T'_{n,i}\subseteq T'_{n,j}$.

Proof.  Clear.

Claim 2. Suppose $T^n$ is an $x_n$-sized tree.  Then $(\forall j\in\omega)\allowbreak(T'_{n,j}$
is an $x_{n,j}$-sized tree).

Proof:  It is clear that $T'_{n,0}$ is an $x_{n,0}$-sized tree.  
Assume that $T'_{n,j}$ is an $x_{n,j}$-sized tree.  Fix $t\in\omega$. 

Case 1: 
 $t<m(n,j)$.

We have that

\medskip

\centerline{ $T'_{n,j+1}\cap {}^t\omega=
T'_{n,j}\cap{}^t\omega$}

\medskip

\noindent and so 

\medskip

\centerline{$\vert T'_{n,j+1}\cap {}^t\omega\vert\leq x_{n,j}(t)\leq
x_{n,j+1}(t)$.}

\medskip

Case 2: $t\geq m(n,j)$.  

We have 

\medskip

\centerline{$T'_{n,j+1}\cap{}^t\omega\subseteq(T'_{n,j}\cap{}^t\omega)\cup (T_{n,j}\cap
{}^t\omega)$.}

\medskip

\noindent  Therefore
we have 

\medskip

\centerline{$\vert T'_{n,j+1}\cap{}^t\omega\vert\leq 2x_{n,j}(t)\leq x_{n,j+1}(t)$.}

\medskip

  The
Claim is established.

For each $n\in\omega$, using Claim 2 and Lemma 4.7 with $k_{n,j}$ here equal to $m(j)$ there,
we my find an increasing sequence of integers $\langle k_{n,j}\,\colon\allowbreak j\in\omega\rangle$ and $T^{n+1}\in V$ such that $k_{n,0}=0$ and
$(\forall j>0)\allowbreak(k_{n,j}>j)$ and
if $T^n$ is an $x_n$-sized tree, then $T^{n+1}$ is
an $x_{n+1}$-sized tree such that
for all $\eta\in{}^{<\omega}\omega$, we have 

\medskip

{\centerline{$(\forall j>0)\allowbreak
(\exists i< j)\allowbreak
(\eta\restr k_{n,j+1}\in T'_{n,k_{n,i}})$ implies $\eta\in T^{n+1}.$}}

\medskip

This completes the construction of $\langle T^n\,\colon\allowbreak n\in\omega\rangle$ and
$\langle T'_{n,j}\,\colon\allowbreak j\in\omega$, $n\in\omega\rangle$.

Applying mathematical induction, we have that each $T^n$ is in fact an $x_n$-sized tree.

Claim 3. $T^n\subseteq T^{n+1}$ for every $n\in\omega$.

Proof: By Claim 1 we have that $T^n\subseteq T'_{n,i}$ for every $i\in\omega$.
By the definition of $T^{n+1}$ we have that

\medskip

\centerline{
$T^{n+1}\supseteq\bigcap\{T'_{n,k_{n,i}}\,\colon\allowbreak i\in\omega\}\supseteq
\bigcap\{T'_{n,i}\,\colon\allowbreak i\in\omega\}\supseteq T^n$.}

\medskip

The Claim is established.

Applying Lemma 4.7 again, with $k_n$ here equal to $m(n)$ there,
we obtain an increasing sequence of integers $\langle k_n\,\colon\allowbreak n\in\omega\rangle$ and a
$z$-sized tree $T^*\in V$ such that $(\forall n\in\omega)\allowbreak
(n<k_n)$ and for every $\eta\in{}^{<\omega}\omega$,
we have that

\medskip

{\centerline{$(\forall n>0)(\exists i< n)(
\eta\restr k_{n+1}\in T^{k_i})$ implies $\eta\in T^*$.}}

\medskip

Notice that $T^0\subseteq\bigcap\{T^n\,\colon\allowbreak n\in\omega\}\subseteq T^*$.

Now we verify that $\langle T^n\,\colon\allowbreak n\in\omega\rangle$ and $T^*$ satisfy the
remaining conclusions of the Lemma.  Accordingly, fix 
$n\in\omega$ and $j\in\omega$ and
$g\in[T_{n,j}]$.
Let 

\medskip

\centerline{$k={\rm max}(k_n,\allowbreak{\rm max}\{k_{n,j'}\,\colon\allowbreak j'\leq j+1\},
\allowbreak{\rm max}\{m(n,j')\,\colon\allowbreak
j'\leq j+1\})$.}

\medskip

  Fix $\eta\in T_{n,j}$ extending $g\restr k$ and assume that
 $\eta\restr k\in T^n\cap T^*$.

Claim 4. $\eta\in T^{n+1}$.

Proof:  It suffices to show 

\medskip

\centerline{$(\forall j'>0)\allowbreak
(\exists i< j')\allowbreak(\eta\restr k_{n,j'+1}\in T'_{n,k_{n,i}})$.}

\medskip

  Fix $j'\in\omega$
and let $i={\rm min}(j,j')$.

Case 1: $j'\leq j$. 

Because $k_{n,j'+1}\leq k$ we have that 
 $\eta\restr k_{n,j'+1}\in T^n\subseteq T'_{n,k_{n,i}}$, as required.

Case 2: $j<j'$.

It suffices to show that $\eta\restr k_{n,j'+1}\in T'_{n,k_{n,j}}$. 
Because $g\restr k=\eta\restr k\in T^n$ and $m(n,j)\leq k$, we have that
 $g\restr m(n,j)\in T^n\subseteq T'_{n,j}$.  Because
we have $\eta\in
T_{n,j}$ and $\eta\restr m(n,j)=g\restr m(n,j)\in T'_{n,j}$, we know
 by the definition of
$T'_{n,j+1}$ and Claim~1 that
$\eta\in T'_{n,j+1}\subseteq T'_{n,k_{n,j}}$.

Claim 4 is established.

Claim 5. $\eta\in T^*$.

Proof: It suffices to show $(\forall i>0)\allowbreak(\exists i'< i)
\allowbreak(\eta\restr k_{i+1}\in T^{k_{i'}})$.  Towards this end, fix $i\in\omega$.

Case 1:  $i\leq n$.

Because $\eta\restr k\in T^*$ and $\eta$ extends $g\restr k$,
 we have $g\restr k_{i+1}\in T^*$ and hence
we may take $i'< i$ such that $g\restr k_{i+1}\in T^{k_{i'}}$. But we also have
 $\eta\restr k_{i+1}=g\restr k_{i+1}$, so we finish Case~1.

Case 2: $n<i$.

We let $i'=i-1$. 
By Claim 4 we have $\eta\restr k_{i+1}\in T^{n+1}$, and by Claim 3 we have that
$T^{n+1}\subseteq T^{k_n}\subseteq T^{k_{i'}}$.

Claim 5 and the Lemma are established.

Using Shelah's terminology, we have by  Lemmas 4.7, 4.8, and 4.9 
 that $(D,R)$ from [13, Definition VI.2.9A] is a smooth
strong covering model [13, Definition VI.1.7(2)] (see [13, Claim VI.1.8(1)].

The proof of the following Theorem is
 [13, proof of Theorem VI.1.12] for the case of [13, Definition VI.2.9A].

\proclaim Theorem 4.10.  Suppose $\langle P_\eta\,\colon\allowbreak\eta\leq\kappa\rangle$ is a countable support iteration based on\/
 $\langle Q_\eta\,\colon\allowbreak\eta<\kappa\rangle$ and suppose\/ {\rm $(\forall\eta<\kappa)\allowbreak({\bf 1}\forces_{P_\eta}``Q_\eta$ 
is proper and has the Sacks property'').} 
Then $P_\kappa$ has the Sacks property.

Proof: The proof proceeds by induction on $\kappa$.
Since no counterexample can first appear when ${\rm cf}(\kappa)$ is uncountable, we may, using standard arguments, 
assume $\kappa$ is either 2 or $\omega$. The case $\kappa = 2$ is easily handled using Lemma 4.6, so assume $\kappa =\omega$.

 Suppose  that 
  $\lambda$ is a sufficiently large regular cardinal and $N$ is a countable elementary
sustructure of $H_\lambda$ and $1\ll z$ and 
${\bf 1}\forces_{P_\omega}``\zeta\in{}^\omega\omega$,'' and
$\{P_\omega,  z\}\in N$ and the $P_\omega$-name $\zeta$ is in $N$ and
$p\in P_\omega\cap N$.

Fix $p'$ and
 $\langle (p_n,\zeta_n)\,\colon\allowbreak n\in\omega\rangle\in N$ as in Lemma 3.1, using $\zeta$ for $f$ and
$\zeta_n$ for $f_n$.

Fix $y\ll z$ such that $1\ll y$ and $y\in N$. 

Let $\Omega=\{x\in N\,\colon\allowbreak 
 1\ll x\ll y\}$. Let $\langle y_n\,\colon\allowbreak n\in\omega\rangle$ enumerate $\Omega$. 
 Build $\langle x^*_n\,\colon\allowbreak n\in\omega\rangle$ as follows.  Fix $x^*_0\in\Omega$, and for each $n\in\omega$ choose $x^*_{n+1}\in \Omega$ such that
$x^*_n\ll x^*_{n+1}$ and $y_{n}\ll x^*_{n+1}$.
Also build $\langle x_n\,\colon\allowbreak n\in\omega\rangle$, $\langle y^*_n\,\colon\allowbreak n\in\omega\rangle$, $\langle x'_n\,\colon\allowbreak n\in\omega\rangle$, and $\langle y'_n\,\colon\allowbreak n\in\omega\rangle$ 
 sequences of elements of $\Omega$ such that 
 for each $n\in\omega$ we have
$x^*_n\ll x_n\ll y^*_n\ll x'_n \ll y'_n\ll x^*_{n+1}$.

For each $n\in\omega$ let $\langle T_{n,j}\,\colon\allowbreak j\in\omega\rangle$ list,
with infinitely many repetitions,  all 
$T'\in N$ such that  there is some $y^*\ll x^*_{n+1}$ such that $T'$ is a
$y^*$-sized tree, and build $\langle x_{n,j}\,\colon\allowbreak j\in\omega\rangle$ a sequence of elements of
$\Omega$ such that for every $j\in\omega$ we have  that
$x_n\ll x_{n,j}\ll x_{n,j+1}\ll x^*_{n+1}$ and $T_{n,j}$ is an $x_{n,j}$-sized tree.

Working in $V[G_{P_\omega}]$, use Lemma 4.9 to choose $T^*$  a $z$-sized tree  and
 $\langle T^n\,\colon\allowbreak n\in\omega\rangle\in V[G_{P_\omega}]$ a sequence  such that
$T^*\in V$ and
$T^0\subseteq T^*$ and $\zeta_0\in [T^0]$ and for every $n\in\omega$ we have that
$T^n\in V$ is an $x_n$-sized tree and $T^n\subseteq T^{n+1}$ and
 for every $j\in\omega$ and every 
$g\in[T_{n,j}]$  there is $k\in\omega$ such that for
every $\eta\in T_{n,j}$ extending $g\restr k$, if
 $\eta\restr k\in T^n\cap T^*$ then
$\eta\in T^{n+1}\cap T^*$. We may assume the $P_\omega$-names $T^*$ and $T^n$ are in $N$ for each $n$.

Note that the reason we worked in $V[G_{P_\omega}]$ rather than in $V$ in the previous paragraph is because 
we wish to allow $g$ to range over $[T_{n,j}]$ with the brackets interpreted in $ V[G_{P_\omega}]$ (i.e., $g$ need not be in
$V$).

Using the induction hypothesis, for every $n\in\omega$, let $F_{n,0}$ and $F_{n,1}$ and $F_{n,2}$ be $P_{n}$-names  such that

\medskip

(A)  
 ${\bf 1}\forces_{P_n}``F_{n,0}$ and $F_{n,1}$ and
$F_{n,2}$ are functions
 each of whose domains
 is equal to $Q_n$, such that

\centerline{$(\forall q'\in Q_n
 )\allowbreak (F_{n,0}(q')$ is an $F_{n,1}(q')$-sized tree}

\centerline{and $F_{n,1}(q')\ll y^*_{n}$ and
$F_{n,2}(q')\leq q'$}

\centerline{and $F_{n,2}(q')\forces`\zeta_{n+1}\in[ F_{n,0}(q')]$').''}

\medskip

 We may assume that the names $F_{n,0}$ and
$F_{n,1}$ and $F_{n,2}$ are in $N$.

For each $n\in\omega$ we may, in $V[G_{P_{n}}]$, use Lemma 4.8 to choose 
$ T_n^\#$ a $x'_n$-sized tree  and
$\langle k^n_i\,\colon\allowbreak i\in\omega\rangle$ an increasing sequence of
integers 
such that 
 $T^n\subseteq T_n^\#$ and for every $\eta\in T^n$ and every $i\in\omega$ and every
$\nu\in F_{n,0}(p_{k^n_i}(n))$,
 if 
${\rm length}(\eta)\geq k^n_i$ and 
$\nu$ extends $\eta$, then $\nu\in T_n^\#$. 

\medskip

Using the induction hypothesis and Lemma 4.6, we may take  $\tilde T_n$ such that
$\tilde T_n\in V$ is $y_n'$-sized tree and $T_n^\#\subseteq \tilde T_n$.

We may assume the $P_{n}$-names $\tilde T_n$ and 
$\langle k^n_i\,\colon\allowbreak i\in\omega\rangle$ are in $N$.

Claim 1.  We may be build $\langle r_n\,\colon\allowbreak n\in\omega\rangle$ such that
 for every $n\in\omega$ we have that the following hold:

(1) $r_n\in P_{n}$ is $N$-generic, and

(2) $r_{n+1}\restr n=r_n$, and

(3) $r_n\forces``\zeta_n\in [T^n]\cap [T^*]$,'' and

(4) $r_n\leq p'\restr n$.''

\medskip

Proof: By induction on $n$.  For $n=0$ we have nothing to prove. Suppose we have $r_n$.

  By (A) and the definition of $\tilde T_n$  we have that

\medskip

(B)
$r_n\forces``T^n\subseteq\tilde T_n$ ''

and 

(C) $r_n\forces``$for every $\eta\in T^n$ and every $i\in\omega$ and every
$\nu\in F_{n,0}(p_{k^n_i}(n))$,

\centerline{if 
${\rm length}(\eta)\geq k^n_i$ and 
$\nu$ extends $\eta$ then $\nu\in\tilde T_n$.''}

\medskip

By (C) and the fact that, by the induction hypothesis, 
we know $r_n\forces``\zeta_n\in[T^n]$,'' we have that

(D) $r_n\forces``(\forall j\in\omega)(F_{n,2}(p_{k_j^n}(n))\forces`(\forall\nu\in
F_{n,0}(p_{k^n_j}(n)))\allowbreak($if $\nu$ extends $\zeta_n\restr k^n_j$ then $\nu\in\tilde T_n))$'\thinspace.''

We have

\medskip

{\centerline{$r_n\forces``\tilde T_n\in N[G_{P_{n}}]$.''}}

\medskip

We also have

\medskip

{\centerline{$r_n\forces``\tilde T_n\in V$.''}}

\medskip

Therefore, because $r_n$ is $N$-generic, we have

\medskip

{\centerline{$r_n\forces``\tilde T_n\in N$.''}}

\medskip

Therefore there is a $P_n$-name $m$ such that

\medskip

{\centerline{$r_n\forces``\tilde T_n=T_{n,m}$ and $m >n$.''}}

\medskip

Using this fact along with the fact that $\langle T^n\,\colon\allowbreak 
n\in\omega\rangle$ and $T^*$ were chosen as in the conclusion of  Lemma 4.9 
and also using the fact that $r_n\forces``\zeta_n\in[T^n]\subseteq[\tilde T_n]$,'' we may choose
 $k$ to be a $P_{n}$-name for an integer such that

\medskip

(E) $r_n\forces``(\forall\eta\in \tilde T_n)\allowbreak($if $\eta$ extends $\zeta_n\restr k$
 and $\eta\restr k\in T^n\cap T^*$ 
then $\eta\in T^{n+1}\cap T^*)$.''

\medskip

Choose $j$  to be a $P_{n}$-name for an integer such that
$r_n\forces``k^n_j\geq k$.''

Subclaim 1. 
$r_n\forces``F_{n,2}(p_{k^n_j}(n))\forces
`\zeta_{n+1}\in [\tilde T_n]$.'\thinspace''

Proof.  It suffices to show

\medskip 

$r_n\forces``F_{n,2}(p_{k^n_j}(n))\forces
`(\forall j'>j)\allowbreak(\zeta_{n+1}\restr k^n_{j'}\in\tilde T_n)$.'\thinspace''

\medskip

Fix $j'$ a $P_{{n+1}}$-name for an integer such that 

\medskip

$r_n\forces``F_{n,2}(p_{k^n_j}(n))\forces` j'>j$.'\thinspace''

\medskip

We know by the induction hypothesis that $r_n\forces``\zeta_n\in[T^n]$.'' Therefore

\medskip

(F) $r_n\forces``\zeta_n\restr k^n_{j}\in T^n$.''

\medskip

By the definition of $\langle p_i\,\colon\allowbreak i\in\omega\rangle$, we have

\medskip

(G) $r_n\forces``p_{k^n_j}(n)\forces`\zeta_n\restr k^n_j=
 \zeta_{n+1}\restr k^n_j$.'\thinspace''

\medskip

By (A) we have

\medskip

(H) $r_n\forces``F_{n,2}(p_{k^n_j}(n))\forces`\zeta_{n+1}\in
[F_{n,0}(p_{k^n_j}(n)])$.'\thinspace''

\medskip

Combining (F), (G), (H), and the definition of $\tilde T_n$, we have that

\medskip

{\centerline{$r_n\forces``F_{n,2}(p_{k^n_j}(n))\forces`\zeta_{n+1}\restr 
k^n_{j'}\in\tilde T_n$.'\thinspace''}}

The Subclaim is established.

Subclaim 2. $r_n\forces``F_{n,2}(p_{k^n_j}(n))\forces`\zeta_{n+1}\in [T^{n+1}]\cap [T^*]$.'\thinspace''

Proof: By (E) we have

\medskip

(I) $r_n\forces``(\forall\eta\in \tilde T_n)\allowbreak($if $\eta$ extends $\zeta_n\restr k^n_j$ and $\eta\restr k^n_j\in T^n\cap T^*$

\noindent then
$\eta\in T^{n+1}\cap T^*$).''

\medskip

Work in $V[G_{P_{n}}]$ with $r_n\in G_{P_{n}}$.
Fix $\eta\in\tilde T_n$ and
 suppose 

\medskip

\centerline{$F_{n,2}(p_{k^n_j}(n))\forces``
\eta$ is an initial segment of $\zeta_{n+1}$ with
${\rm lh}(\eta)\geq k^n_j$.''}

\medskip

To establish the Subclaim, it suffices to show

\medskip

(J) $F_{n,2}(p_{k^n_j}(n))\forces``\eta\in T^{n+1}\cap T^*$.''

\medskip

By the definition of $\langle p_i\,\colon\allowbreak i\in\omega\rangle$ we have

\medskip

{\centerline{$p_{k^n_j}(n)\forces``\eta\restr k^n_j= \zeta_{n+1}\restr k^n_j=\zeta_n\restr k^n_j$.''}}

\medskip

Hence by the fact that Claim 1 holds for the integer $n$ we have

\medskip

(K) $p_{k^n_j}(n)\forces``\eta\restr k^n_j\in T^n\cap T^*$.''

\medskip

By Subclaim 1, (I), (K), and the fact that
$F_{n,2}(p_{k^n_j}(n))\leq
p_{k^n_j}(n)$ we obtain (J).

Subclaim 2 is established.

To complete the induction establishing Claim 1, we  take
$r_{n+1}\in P_{{n+1}}$ such that $r_{n+1}\restr n=r_n$ and
$r_{n+1}$ is $N$-generic and $r_n\forces``r_{n+1}(n)\leq
F_{n,2}(p_{k^n_j}(n))$.''

Claim 1 is established.

Define $q'$ by

\medskip

\centerline{$q'=\bigcup\{r_n\restr n\,\colon\allowbreak n\in\omega\}$.}

\medskip

 By Claim 1 we have
that

\medskip

\centerline{$q'\forces``$for every $n\in\omega$ we have 
$\zeta_n\in[T^*]$ and $\zeta_n\restr n=\zeta\restr n$,}

\centerline{and therefore $\zeta\in[T^*]$.''}

\medskip

The Theorem is established.

\section{The Laver Property}

In this section, we present Shelah's proof that the Laver property is preserved by countable support iteration of proper forcing.

\proclaim Definition 5.1. Suppose $f\in {}^\omega(\omega-\{0\})$ and\/ ${\bf 1}\ll f$.
We say that\/ $T$ is an $f$-tree iff $T$ is a tree and\/ $(\forall\eta\in T)\allowbreak(\forall n\in{\rm dom}(\eta))\allowbreak(\eta(n)<f(n))$.

\proclaim Definition 5.2.
We say that\/ $P$ is $f$-preserving iff whenever $z$ is
in ${}^\omega(\omega-\{0\})$ and\/ $1\ll z$ then 

${\bf 1}\forces_P``(\forall g\in{}^\omega\omega)(g\leq f$ implies there exists $H\in V$ such that $H$ is
a $z$-sized $f$-tree  and $g\in[ H])$.''

\medskip

\proclaim Definition 5.3. We say that\/ $P$ has the Laver property iff for every $f\in{}^\omega(\omega-\{0\})$ such that\/ ${\bf 1}\ll f$ we have that\/ $P$ is $f$-preserving.

The following Theorem is [13, Claim VI.2.10C(2)].

\proclaim Theorem 5.4. $P$ has the Sacks property iff\/ $P$ has the Laver property and\/ $P$ is ${}^\omega\omega$-bounding.

Proof:  We first assume that $P$ has the Sacks property and we show that $P$ is ${}^\omega\omega$-bounding.  Given $p\in P$ and a name $f$ such that $p\forces``f\in{}^\omega\omega$,'' take $z\in{}^\omega(\omega-\{0\})$ such that ${\bf 1}\ll z$ and use the fact that $P$ has the Sacks property to obtain $q\leq p$ and a $z$-sized tree $H$ such that
$q\forces``f\in[H]$.'' For every $n\in\omega$ let 

\medskip

\centerline{$g(n)={\rm max}\{\eta(n)\,\colon\allowbreak\eta\in H$ and ${\rm lh}(\eta)>n\}$.}

\medskip

 Then we have $q\forces``f\leq g$.'' This establishes the fact that $P$ is ${}^\omega\omega$-bounding.

It is clear that if $P$ has the Sacks property, then it has the Laver property.

Finally we assume that $P$ has the Laver property and is ${}^\omega\omega$-bounding, and we show that $P$ has the Sacks property.
So suppose that $p\in P$ and ${\bf 1}\ll z$ and $p\forces``g\in{}^\omega\omega$.''
Using the fact that $P$ is ${}^\omega\omega$-bounding, take $p'\leq p$ and  $f\in{}^\omega\omega$ such that
$p'\forces``g\leq f$.'' Using the fact that $P$ has the Laver property, take $q\leq p'$ and $H$ a
$P$-name such that

\medskip

\centerline{$q\forces``H$ is a $z$-sized $f$-tree and $g\in[H]$ and $H\in V$.''}

\medskip

The Theorem is established.

The following is [13, Conclusion VI.2.10D].

\proclaim Theorem 5.5. Suppose $\langle P_\eta\,\colon\allowbreak\eta\leq\kappa\rangle$ is a countable support iteration based on\/
 $\langle Q_\eta\,\colon\allowbreak\eta<\kappa\rangle$ and suppose\/ {\rm $(\forall\eta<\kappa)\allowbreak({\bf 1}\forces_{P_\eta}``Q_\eta$ 
is proper and has the Laver property.'')}  Then $P_\kappa$ has the Laver property.

Proof: Fix $f\in{}^\omega\omega$ such that ${\bf 1}\ll f$.  Repeat the proofs of Lemma 4.6 through Theorem 4.10 with ``tree'' replaced by ``$f$-tree.'' The Theorem is established.

\section{$(f,g)$-bounding}

In this section we establish the preservation of $(f,g)$-bounding forcing.  For an exact formulation, see
 Theorem 6.5 below.  This proof is due to Shelah, of course; see [13, Conclusion VI.2.11F].

The following Definition corresponds to [13, Definition VI.2.11A].

\proclaim Definition 6.1. We say that $T$ is an $(f,g)$-corseted tree iff

(0) $T\subseteq{}^{<\omega}\omega$ is a tree with no terminal nodes, and

(1)  $f$ and $g$ 
are functions with domain $\omega$, and

(2) $(\forall n\in\omega)\allowbreak(f(n)\in\{r\in{\bf R}\,\colon
\allowbreak 1<r\}\cup\{\omega\})$, and

(3) $(\forall n\in\omega)\allowbreak(g(n)\in\{r\in {\bf R}\,\colon\allowbreak
1<r\}\cup\{\aleph_0\})$, and

(4) $f$ and $g$ diverge to infinity, and

(5) $(\forall\eta\in T)(\forall i\in{\rm dom}(\eta))(\eta(i)<f(i))$, and

(6) $(\forall n\in\omega)(\vert\{\eta(n)\,\colon\eta\in T$ and $n\in{\rm dom}(\eta)\}\vert\leq g(n))$.

\proclaim Definition 6.2. Suppose that $f$ and $g$ are functions as in Definition 6.1.
We say that $P$ is $(f,g)$-bounding iff\/ {\rm ${\bf 1}\forces_P``(\forall h\in {}^\omega\omega)\allowbreak [(\forall n\in\omega)\allowbreak(h(n)<f(n))$ implies $(\exists T\in V)\allowbreak
(T $ is an $(f,g)$-corseted tree and $h\in[T])]$.''}

The following Lemma is the analogue  of Lemma 4.6.

\proclaim Lemma 6.3.  $P$ is $(f^{g^k},g^{1/k})$-bounding for infinitely many $k\in\omega$ iff
whenever $x<z$ are positive rational numbers and\/ $\gamma\in\omega$ then\/ {\rm ${\bf 1}\forces``$
if $T $ is an 
$(f^{g^\gamma}, g^x)$-corseted tree then $(\exists H\in V)\allowbreak
(H$ is an $(f^{g^\gamma}, g^z)$-corseted tree and $T\subseteq H)$.''}

Proof.  We prove the non-trivial direction.
Fix an integer $k$ such that $k> x$  and $P$ is $(f^{g^{\gamma+k}},g^{1/{(\gamma+k)}})$-bounding and $k>1/(z-x)$.
Let $X=\{n\in\omega\,\colon\allowbreak g(n)=\aleph_0\}$.

For every $m\in\omega-X$ define 

\medskip

\centerline{${\cal T}_m=\{S\subseteq\omega\,\colon\allowbreak
{\rm sup}(S)\leq f(m)^{g(m)^\gamma}$ and $\vert S\vert<g(m)^x\}$.}

\medskip

For every $m\in\omega-X$ define

\medskip

\centerline{${\cal T}'_m=\{i\in\omega\,\colon\allowbreak
 i\leq (f(m)^{g(m)^\gamma}+1)^{g(m)^k}\}$.}

\medskip

Because $x<k$ we may choose, for each integer $m$ not in $X$, a one-to-one mapping $h_m$
 from ${\cal T}_m$
into ${\cal T}'_m$.

Define

\medskip

\centerline{${\cal T}=\{\xi\in{}^{<\omega}\omega\,\colon(\forall m\in \omega-X)(\xi(m)\in{\cal T}'_m)$}

\centerline{and $(\forall m\in X)(\xi(m)=1)\}$.}

\medskip

In $V[G_P]$ let $\zeta\in[{\cal T}]$ denote the function defined by

\medskip

\centerline{$(\forall m\in\omega-X)\allowbreak(\zeta(m)=h_m(\{\eta(m)\,\colon\allowbreak\eta\in T$ and $m\in{\rm dom}(\eta)\}))$}

\centerline{and $(\forall m\in X)\allowbreak(\zeta(m)=1)$.}

\medskip

Because $P$ is
$(f^{g^{\gamma+k}},g^{1/{(\gamma+k)}})$-bounding, 
we may take $H'\in V$ such that $H'$ is an $(f^{g^{\gamma+k}},g^{1/{(\gamma+k)}})$-corseted tree and
$\zeta\in[H']$. Define $H$ by

\medskip

\centerline{$H(m)=\bigcup\{h^{-1}_m(t)\,\colon\allowbreak(\exists\eta\in H')( t=\eta(m))$ and $t\in {\rm range}(h_m)\}$ for $m\in\omega-X$,}

\centerline{and $H(m)=\omega$ for $m\in X$.}

\medskip

When $g(m)$ is finite, we have

\medskip

{\centerline{$\vert H(m)\vert\leq\vert H'(m)\vert\cdot{\rm max}\{\vert h_m^{-1}(t)\,\colon\allowbreak t\in{\rm range}(h_m)\vert\}$}}

{\centerline{$\leq
 g^x(m)\cdot g^{1/{(\gamma+k)}}(m)< g^z(m)$.}}

\medskip

We have that $H$ is an $(f^{g^\gamma},g^z)$-corseted tree and ${\bf 1}\forces``T\subseteq H$.''  The Lemma is established.

The following Lemma is the analogue of Lemma 4.7.

\proclaim Lemma 6.4.  Suppose\/  $\langle r_n\,\colon\allowbreak n\in\omega\rangle$ is a
bounded sequence of positive rational numbers and\/ $y\in{\bf Q}$ and\/ ${\rm sup}\{r_n\,\colon\allowbreak
n\in\omega\}<y$. Suppose $P$ is a forcing such that\/ {\rm
$V[G_P]\models``$for every countable $X\subseteq V$ there is a countable $Y\in V$ such that
$X\subseteq Y$.''}.  Suppose in $V[G_P]$ we have\/ {\rm $(\forall n\in\omega)\allowbreak(T_n\in V$ is an 
$(f,g^{r_n})$-corseted tree).}  
Then in  $V[G_P]$ there is an $(f,g^{y})$-corseted tree $T^*\in V$ and an increasing sequence of integers 
$\langle k_n\,\colon\allowbreak n\in\omega\rangle$ such that\/ $k_0=0$ and\/
$(\forall i>0)\allowbreak(i<k_i)$ and
for every $\eta\in{}^{<\omega}\omega$ we have 

\medskip

\centerline{$(\forall t\in{\rm dom}(\eta))(\exists j\in\omega)(\exists\nu\in T_{k_j})(k_j\leq t$ and $\nu(t)=\eta(t))$}

\centerline{implies $\eta\in T^*).$}

\medskip

Proof: The proof is similar to the proof of Lemma 4.7.
We note the following modifications. We must choose $x\in{\bf Q}$ such that ${\rm sup}\{r_n\,\colon\allowbreak n\in\omega\}<x<y$. 
By recursion choose $\langle k_n\,\colon\allowbreak n\in\omega\rangle$ an increasing sequence of integers such that $k_0=0$ and
$(\forall n>0)\allowbreak(\exists j\in\omega)\allowbreak(k_n\leq j$ implies
$(n+1)g(j)^{x}\leq g(j)^{y})$.

The definition of $T^*$ is changed to $T^*=\{\eta\in{}^{<\omega}\omega\,\colon(\forall t\in{\rm dom}(\eta))(\exists j\in\omega)(\exists\nu\in S'_{k_j})(k_j\leq t$ and $\nu(t)=\eta(t))\}$.

Clearly $T^*$ is a tree.

Claim. $T^*$ is an $(f,g^y)$-corseted tree.

Proof: Fix $t\in\omega$. 

Case 1:  $t\geq k_1$.

Choose $m\in\omega$ such that $k_m\leq t<k_{m+1}$.  We have that

\medskip

\centerline{
$\vert \{\eta(t)\,\colon\eta\in T^*$ and $t\in{\rm dom}(\eta)\}\vert=
\Sigma_{j\leq m}\vert \{\eta(t)\,\colon\allowbreak\eta\in T_{k_j}$ and
$t\in{\rm dom}(\eta)\}\vert$}

\centerline{$\leq (m+1)g^{x}(t)\leq g^{y}(t)$.}

\medskip

Case 2: $t<k_1$.

We have $\{\eta(t)\,\colon \eta\in T^*\}=\{\eta(t)\,\colon \eta\in T\}$, so it follows that
$\vert H(t)\vert\leq g^y(t)$.

The Claim is established.

The other requirements of the Lemma are the same as in the proof of Lemma 4.7. The Lemma is established.

The following is [13, Conclusion VI.2.11F].

\proclaim Theorem 6.5.   Suppose $\langle P_\eta\,\colon\allowbreak\eta\leq\kappa\rangle$ is a countable support iteration
based on $\langle Q_\eta\,\colon\allowbreak\eta<\kappa\rangle$ and suppose 
that for every $\eta<\kappa$ we have that\/ {\rm
${\bf 1}\forces``$for infinitely many $k\in\omega$ we have that $Q_\eta$ is proper and
$(f^{g^k},g^{1/k})$-bounding.''} Then $P_\kappa$ is
$(f^{g^k},g^{1/k})$-bounding for every positive $k\in\omega$.

Proof: The same as Theorem 4.10, with ${}^\omega(\omega-\{0\})$ replaced by ${\bf Q}$,
and $\ll$ replaced by $<$, and $x$-sized tree replaced by $(f^{g^\gamma},g^x)$-corseted tree.

\section{$P$-point property}

In this section we define the $P$-point property and prove that it is preserved by countable support iteration of proper
forcings.  This is due to Shelah [13, Conclusion VI.2.12G].

\proclaim Definition 7.1. Suppose $n\in\omega$ and $x\in{}^\omega(\omega-\{0\})$ is strictly increasing.
We say that $(j,k,m)$ is an $x$-bound system above $n$ iff each of the following holds:

(0) $k\in\omega$, and

(1)   $j$ and $m$ are functions from $k+1$ into $\omega$, and

(2) $j(0)>x(n+m(0)+1)$, and

(3) $(\forall l< k)\allowbreak(j(l+1)>x(j(l)+m(l+1)+1))$.

\proclaim Definition 7.2. Suppose $n\in\omega$ and $x\in{}^\omega(\omega-\{0\})$ is strictly increasing
and $(j,k,m)$ is an $x$-bound system above $n$ and $T$ is a tree.  We say that $T$ is a
$(j,k,m,\eta)$-squeezed tree iff $T$ has no terminal nodes and
each of the following holds:

(1) ${\rm dom}(\eta)=\{(l,t)\in{}^2\omega\,\colon\allowbreak
l\leq k$ and $t\leq m(l)\}$, and

(2) $(\forall (l,t)\in{\rm dom}(\eta))\allowbreak
(\eta(l,t)\in {}^{j(l)}\omega)$, and

(3) $(\forall\nu\in T)\allowbreak(\exists (l,t)\in{\rm dom}
(\eta))\allowbreak(\nu$ is comparable with $\eta(l,t))$.

It is easy to see that the following Definition is equivalent to [13, Definition VI.2.12A].

\proclaim Definition 7.3. We say that\/ $T$ is $x$-squeezed iff for every $n\in\omega$ there is some $x$-bound system
$(j,m,k)$ above $n$ such that\/ $T$ is $(j,k,m,\eta)$-squeezed for some $\eta$.

In other words, $T$ is $x$-squeezed when, living above any given level of $T$, say $\{\xi\in T\,\colon\allowbreak
{\rm lh}(\xi)=n+1\}$, there is a maximal antichain ${\cal A}$ of $T$ that can be decomposed
as ${\cal A}=\bigcup\{{\cal A}_l\,\colon\allowbreak l\leq k\}$ where each
${\cal A}_l$ is a subset of $\{\xi\in T\,\colon\allowbreak{\rm lh}(\xi)=j(l)\}$ of cardinality at most $m(l)+1$, such that
the levels of ${\cal A}$ are stratified so sparsely that conditions (2) and (3) of Definition 7.1 hold. 
 Notice that for any given $l\leq k$ we may have that
$\{\eta(l,t)\,\colon\allowbreak t\leq m(l)\}$ is a proper superset of ${\cal A}_l$; indeed, it need not even be a subset of $T$.
We could modify Definition 7.2 to require this, but there is no need to do so.

\proclaim Lemma 7.4. Suppose $1\ll x\ll y$ and both $x$ and $y$ are strictly increasing
and\/ $T$ is a $y$-squeezed tree. Then\/ $T$ is an $x$-squeezed tree.

Proof: Every $y$-bound system is an $x$-bound system.

\proclaim Definition 7.5.  We say that $P$ has the $P$-point property iff for every $x\in{}^\omega(\omega-\{0\})$ strictly increasing, we have

{\centerline{${\bf 1}\forces``(\forall f\in{}^\omega\omega)(\exists H\in V)(f\in [H]$ and $H$ is an $x$-squeezed tree).''}}

\medskip

\proclaim Lemma 7.6.  $P$ has the $P$-point property iff for every $x\in
{}^\omega(\omega-\{0\})$  strictly increasing and every
 $p\in P$, if\/ {\rm $p\forces_P``f\in{}^\omega\omega$''}  there are $q\leq p$ and an 
 $x$-squeezed tree
 $H$ such that\/
$q\forces``f\in[H]$.''

Proof: Assume that $P$ has the $P$-point property.  Given $x$, $p$, and $f$, there is $q\leq p$ and $H\subseteq{}^{<\omega}\omega$ such that
$q\forces``f\in[H]$ and $H$ is an $x$-squeezed tree.'' By the Shoenfield Absoluteness Theorem we have that
$H$ is an $x$-squeezed tree.

The other direction is immediate, and so the Lemma is established.

\proclaim Lemma 7.7. Suppose $T$ is an $x$-squeezed tree and $n\in\omega$. Then $T\cap{}^n\omega$ is finite.

Proof. Fix  $(j,k,m)$ an $x$-bound system above $n$ and fix $\eta$ such that
$T$ is a $(j,k,m,\eta)$-squeezed tree.
We have $T\cap{}^n\omega\subseteq\{\eta(s,t)\restr n\,\colon\allowbreak
t\leq j(k)$ and $s\leq m(t)\}$.

The following Lemma is [13, Claim VI.2.12B(1)].

\proclaim Lemma 7.8. Suppose that $P$ has the $P$-point
property. Then $P$ is ${}^\omega\omega$-bounding.

Proof: Suppose $p\in P$ and $p\forces``f\in{}^\omega\omega$.'' Pick $x\in{}^\omega
(\omega-\{0\})$ such that $1\ll x$, and take $q\leq p$ and $H$ an $x$-squeezed
tree such that 
$q\forces``f\in[H]$.''  By Lemma 7.7 we may define $h\in{}^\omega\omega$ by
$(\forall n\in\omega)\allowbreak(h(n)={\rm max}\{\nu(n)\,\colon\allowbreak
\nu\in H$ and $n\in{\rm dom}(\nu)\})$. Clearly $q\forces``f\leq h$,'' and
the Lemma is established.

The following Lemma is [13, Claim VI.2.12B(2)].

\proclaim Lemma 7.9. Suppose $P$ has the Sacks property. 
Then $P$ has the $P$-point property.

Proof. Suppose $x\in{}^\omega(\omega-\{0\})$ is strictly increasing and $p\in P$ and
$p\forces``f\in{}^\omega\omega$.'' 
Choose $y\in{}^\omega(\omega-\{0\})$ monotonically non-decreasing such that for $n> x(3)$ 
we have that $y(n)$ is the greatest $t\in\omega$ such that $x(3t)<n$.
Using the Sacks property, choose $q\leq p$ and $H$ a $y$-sized tree such that $q\forces``f\in[H]$.''

Notice that for all $t>0$  we have $y(x(3t))$ is less than or equal to
 the greatest integer $k$
satisfying $x(3k)<x(3t)$, and therefore we have

\begin{itemize}

\item[$(*)$] $(\forall t\in\omega)(y(x(3t))<t)$.

\end{itemize}

Suppose $n>x(3)$.
Let $j$ be such that ${\rm dom}(j)=\{0\}$ and $j(0)=x(2n)+1$ and
let $k=0$; and let $m$ be such that ${\rm dom}(m)=\{0\}$ and $m(0)=\vert H\cap
{}^{j(0)}\omega\vert$. 

Claim: $(j,k,m)$ is an $x$-bound system above $n$.

Proof: We have $x(n+m(0)+1)\leq x(n+1+y(x(2n)+1))\leq x(n+1+y(x(3n)))\leq x(n+1+n-1)<x(2n)+1=j(0)$. 
The first inequality is because $m(0)=\vert H\cap{}^{j(0)}\omega\vert\leq y(j(0))=y(x(2n)+1)$.
The second inequality is because $x$ is strictly increasing and $y$ is monotonically non-decreasing.  The third
inequality is by $(*)$.

The Claim is established.

Define $\eta$ with domain equal to $\{(0.i)\,\colon\allowbreak
i<m(0)\}$ and such that $\langle \eta(0,i)\,\colon\allowbreak i<m(0)\rangle$
enumerates $H\cap{}^{j(0)}\omega$.  Clearly $H$ is a $(j,k,m,\eta)$-squeezed tree, so the Lemma is established.

\proclaim Lemma 7.10. Suppose $y\in{}^\omega(\omega-\{0\})$ is strictly increasing and $T$ and $T'$ are $y$-squeezed trees.
Then $T\cup T'$ is a $y$-squeezed tree.

Proof:  Given $n\in\omega$, choose $(j,k,m,\eta)$  such that
$(j,k,m)$ is a $y$-bound systems above $n$ and
$T$ is $(j,k,m,\eta)$-squeezed. Let $h=j(k)$. Choose $(j',m',k')$ a $y$-bound system above $h$ and choose $\eta'$
such that $T'$ is $(j',\allowbreak k' ,\allowbreak m',\eta')$-squeezed.
We proceed to fuse $(j,k,m,\eta)$ with $(j',\allowbreak k'\allowbreak m',\allowbreak
\eta')$. 
For every 
$l\leq k$ let $j^*(l)=j(l)$ and for every $l$ such that $k<l\leq k+k'+1$ let
$j^*(l)=j'(l-k_n-1)$.
Let $k^*=k+k'+1$.
For every 
$l\leq k$ let $m^*(l)=m(l)$ and for every $l$ such that $k<l\leq k+k'+1$ let
$m^*(l)=m'(l-k-1)$.
For every 
$l\leq k$ and $\beta\leq m(l)$ let $\eta^*(l,\beta)=\eta(l,\beta)$ and for every $l$ such that $k<l\leq k+k'+1$ and
every $\beta\leq m'(l-k-1)$ let
$\eta^*(l,\beta)=\eta'(l-k-1,\beta)$.
It is straightforward to verify that $(j^*,k^*,m^*)$ is a $y$-bound system above $n$ and
that $T\cup T'$ is $(j^*,k^*,m^*,\eta^*)$-squeezed.

The Lemma is established.

\proclaim Definition 7.11. Suppose  $n\in\omega$ and\/
$h\in{}^\omega\omega$ and\/  $y\in{}^\omega(\omega-\{0\})$ is strictly increasing.
Suppose $(j,m,k)$ is a $y$-bound system above $n$.  We say that\/ $(j,m,k)$ is $h$-tight iff\/  $j(k)<h(n)$.
For $T$  a $y$-squeezed tree, we say that $T$ is $h$-tight iff for every $n\in\omega$ there is an
$h$-tight $y$-bound system
$(j,m,k)$  above $n$ such that 
 for some $\eta$ we have that\/ $T$ is $(j,m,k,\eta)$-squeezed for some $\eta$.

The following Lemma should be compared with Lemma 4.6.  Notice the fact that what we prove here is stronger in that the same $y\in{}^\omega(\omega-\{0\})$ is used in both the hypothesis and the conclusion.  This strengthening  is possible by the use of Lemma 7.10.
The proof of Lemma 7.12 is [13, proof of Claim VI.2.12C$(\varepsilon)^+$], except that we have $x=y$.
Thus we are proving a stronger statement than [13, Claim VI.2.12C$(\varepsilon)^+$],
but in fact Shelah likewise proves this stronger statement without saying so.

\proclaim Lemma 7.12.  $P$ has the $P$-point property iff whenever\/ $y\in{}^\omega(\omega-\{0\})$ 
is strictly increasing then\/ {\rm ${\bf 1}\forces``$whenever $T$ is a $y$-squeezed tree
then $(\exists H\in V)\allowbreak(H$ is a $y$-squeezed tree and
$T\subseteq H$).''}

Proof. We prove the non-trivial direction.
Suppose $p\in P$ and $p\forces``T$ is a $y$-squeezed tree.'' Fix $q'\leq p$.
 By Lemma 7.8 we may choose $h\in{}^\omega\omega$ and 
$q\leq q'$ such that 
$q\forces``T$ is $h$-tight.''

Define $z\in{}^\omega(\omega-\{0\})$ by $z(0)=0$ and

\medskip

{\centerline{$(\forall n\in\omega)\allowbreak(z(n+1)= h(z(n))$.}}

\medskip

For every $n\in\omega$ let ${\cal T}_n=\{t\subseteq{}^{< h(n)}\omega\,\colon\allowbreak 
t=T\cap{}^{< h(n)}\omega$ for some $h$-tight $y$-squeezed tree $T\}$. 

Let ${\cal T}=\bigcup\{{\cal T}_n\,\colon n\in\omega\}$. We implicitly fix an isomorphism
from ${}^{<\omega}\omega$ onto ${\cal T}$.

Using the fact that $P$ satisfies the $P$-point property, fix $q^*\leq q$ and ${\cal C}\subseteq{\cal T}$ such that ${\cal C}$ is a
$z$-squeezed tree and $q^*\forces``(\forall n \in\omega)\allowbreak(T\cap{}^{< h(n)}\omega\in{\cal C})$.''

Define $H^*=\bigcup{\cal C}$ and let $H=\{\nu\in H^*\,\colon\allowbreak(\forall n\in\omega)\allowbreak(\exists\eta\in{}^n\omega\cap H^*)
\allowbreak(\eta$ is comparable with $\nu)\}$.

Pick a $z$-bound system $(j^*,k^*,m^*)$ above $n$ and $\eta^*$ such that  ${\cal C}$ is a $(j^*,k^*,\allowbreak m^*,\allowbreak\eta^*)$-squeezed tree.

Fix $n\in\omega$. We show that there is a $y$-bound system $(j,m,k)$ above $n$ such that for some $\eta$ we have
that $H$ is $(j,m,k,\eta)$-squeezed.

Claim 1. For every $\beta\leq m^*(0)$ we have ${\rm ht}(\eta^*(0,\beta))\geq
h(z(n+\beta+1))$.  For every non-zero $\alpha\leq k^*$ and every $\beta\leq m^*(\alpha)$ we have ${\rm ht}(\eta^*(\alpha,\beta))\geq
h(z(j^*(\alpha-1)+\beta+1))$.

Proof:  For every $\beta\leq m^*(0)$ we have ${\rm ht}(\eta^*(0,\beta))=h({\rm rk}_{\cal T}(\eta^*(0,\beta)))=h(j^*(0))\geq
h(z(n+m^*(0)+1))\geq
h(z(n+\beta+1))$.  For every non-zero $\alpha\leq k^*$ and every 
$\beta\leq m^*(\alpha)$ we have ${\rm ht}(\eta^*(\alpha,\beta))=h({\rm rk}_{\cal T}(\eta^*(\alpha,\beta)))=h(j^*(\alpha))\geq
h(z(j^*(\alpha-1)+m^*(\alpha)+1))\geq
h(z(j^*(\alpha-1)+\beta+1))$.

By Claim 1 we may construct $y$-bound systems as follows.
For every  $\beta\leq m^*(0)$, fix an $h$-tight $y$-bound system $(j^{0,\beta},m^{0,\beta},\allowbreak
k^{0,\beta})$ above $z(n+\beta+1)$ along with $\eta^{0.\beta}$ such that
for some $(j^{0,\beta},m^{0,\beta},k^{0,\beta},\eta^{0,\beta})$-squeezed tree $T$ we have $\eta^*(0,\beta)=
{}^{<h(z(n+\beta+1))}\omega\cap T$, and for every non-zero $\alpha\leq k^*$ and $\beta\leq m^*(\alpha)$, 
fix an $h$-tight $y$-bound system $(j^{\alpha,\beta},m^{\alpha,\beta},\allowbreak
k^{\alpha,\beta})$ above $z(j^*(\alpha-1)+\beta+1)$ along with $\eta^{\alpha.\beta}$ such that
for some $(j^{\alpha,\beta},m^{\alpha,\beta},k^{\alpha,\beta},\eta^{\alpha,\beta})$-squeezed tree $T$ we have $\eta^*(\alpha,\beta)=
{}^{<h(z(j^*(\alpha-1)+\beta+1))}\omega\cap T$.

We define

\medskip

\centerline{$\hat\jmath(\alpha,\beta,\gamma)=j^{(\alpha,\beta)}(\gamma)$}

\medskip

\noindent and 

\medskip

\centerline{$\hat k(\alpha,\beta)=k^{(\alpha,\beta)}$}

\medskip

\noindent and 

\medskip

\centerline{$\hat m(\alpha,\beta,\gamma)=m^{(\alpha,\beta)}(\gamma)$}

\medskip

\noindent and for $t\leq\hat m(\alpha,\beta,\gamma)$ let 

\medskip

\centerline{$\hat\eta(\alpha,\beta,\gamma,t)=\eta^{(\alpha,\beta)}(\gamma,t)$.}

\medskip

Claim 2. Suppose $n\in\omega$. Then we have the following:

(1) $\hat\jmath(0,0,0)>y(n+\hat m(0,0,0)+1)$, and

(2) For every $\alpha\leq k^*$ and $\beta\leq m^*(\alpha)$ and $\gamma<\hat k(\alpha,\beta)$ we have
$\hat\jmath(\alpha,\beta,\gamma+1)>y(\hat\jmath(\alpha,\beta,\gamma)+\hat m(\alpha,\beta,\gamma+1)+1)$, and

(3) For every $\alpha\leq k^*$ and every $\beta< m^*(\alpha)$  we have
$\hat\jmath(\alpha,\beta+1,0)>y(\hat\jmath(\alpha,\beta,\allowbreak\hat k(\alpha,\beta))+\hat m(\alpha,\beta+1,0)+1)$, and

(4) For every $\alpha< k^*$ we have
$\hat\jmath(\alpha+1,0,0)>y(\hat\jmath(\alpha,m^*(\alpha),\hat k(\alpha,m^*(\alpha)))+\hat m(\alpha+1,0,0)+1)$.

Proof: Clause (1) holds because $j^{(0,0)}(0)>y(n+m^{(0,0)}(0)+1)$.

Clause (2) holds because $j^{(,\alpha,\beta)}(\gamma+1)>y(j^{(\alpha,\beta)}(\gamma)+m^{(\alpha,\beta)}(\gamma+1)+1)$.

We verify clause (3) as follows.

Case A: $\alpha=0$.

Notice that $j^{0,\beta)}(k^{(0,\beta)})<h(z(n+\beta+1))$ becuase the system $(j^{(0,\beta)},\allowbreak
m^{(0,\beta)},\allowbreak k^{(0,\beta)})$ is $h$-tight above $z(n+\beta+1)$. Notice also that
$j^{(0,\beta+1)}(0)>y(z(n+\beta+2)+m^{(0,\beta+1)}(0)+1)$ because the system $(j^{(0,\beta+1)},\allowbreak
m^{(0,\beta+1)},\allowbreak k^{(0,\beta+1)})$ is above $z(n+\beta+2)$.
 Hence we have
$\hat\jmath(0,\beta+1,0)=j^{(0,\beta+1)}(0)>y(z(n+\beta+2)+m^{(0,\beta+1)}(0)+1)\geq
y(h(z(n+\beta+1))+m^{(0,\beta+1)}(0)+1)\geq
y(j^{(0,\beta)}(k^{(0,\beta)})+m^{(0,\beta+1)}(0)+1)=y(\hat\jmath(0,\beta,\hat k(0,\beta))+\hat m(0,\beta+1,0)+1)$.

Case B: $\alpha>0$.

Notice that $j^{\alpha,\beta)}(k^{(\alpha,\beta)})<h(z(j^*(\alpha-1)+\beta+1))$ becuase the system $(j^{(\alpha,\beta)},\allowbreak
m^{(\alpha,\beta)},\allowbreak k^{(\alpha,\beta)})$ is $h$-tight above $z(j^*(\alpha-1)+\beta+1)$. Notice also that
$j^{(\alpha,\beta+1)}(0)>y(z(j^*(\alpha-1)+\beta+2)+m^{(\alpha,\beta+1)}(0)+1)$ because the system $(j^{(\alpha,\beta+1)},\allowbreak
m^{(\alpha,\beta+1)},\allowbreak k^{(\alpha,\beta+1)})$ is above $z(j^*(\alpha-1)+\beta+2)$.
 Hence we have
$\hat\jmath(\alpha,\beta+1,0)=j^{(\alpha,\beta+1)}(0)>y(z(j^*(\alpha-1)+\beta+2)+m^{(\alpha,\beta+1)}(0)+1)\geq
y(h(z(j^*(\alpha-1)+\beta+1))+m^{(\alpha,\beta+1)}(0)+1)\geq
y(j^{(\alpha,\beta)}(k^{(\alpha,\beta)})+m^{(\alpha,\beta+1)}(0)+1)=y(\hat\jmath(\alpha,\beta,\hat k(\alpha,\beta))+\hat m(\alpha,\beta+1,0)+1)$.

To see that clause (4) holds, we have $\hat\jmath(\alpha+1,0,0)=j^{(\alpha+1,0)}(0)\geq
y(z(j^*(\alpha)+1)+m^{(\alpha+1,0)}(0)+1)\geq y(h(z(j^*(\alpha)))+m^{(\alpha+1,0)}(0)+1)
\geq  y(h(j^*(\alpha))+m^{(\alpha+1,0)}(0)+1)
\geq y(h(z(j^*(\alpha-1))+m^*(\alpha)+1))+m^{(\alpha+1,0)}(0)+1)\geq
y(j^{(\alpha,m^*(\alpha))}(k^{(\alpha,m^*(\alpha))})+m^{(\alpha+1,0)}(0)+1)$.

The first inequality is because the system $(j^{(\alpha+1,0)},m^{(\alpha+1,0)},k^{(\alpha+1,0)})$ is above
$z(j^*(\alpha)+1)$ whence by clause (2) of Definition 7.1 we have the first inequality. 
The second inequality is by definition of the function $z$.  The third inequality is by the fact that $z$ is an increasing function. 
The fourth inequality is because $(j^*,m^*,k^*)$ satisfies clause (2) of Definition 7.1. The fifth inequallity
is because the system $(j^{(\alpha,m^*(\alpha))},m^{(\alpha,m^*(\alpha))},k^{(\alpha,m^*(\alpha))})$ is $h$-tight 
above $z(j^*(\alpha-1)+m^*(\alpha)+1)$.

The Claim is established.

Claim 3.  Suppose $n\in\omega$ and $\nu\in H$. Then there are $\alpha\leq k^*$ and $\beta\leq m^*(\alpha)$ and
$\gamma\leq\hat k(\alpha,\beta)$ and $\delta\leq \hat m(\alpha,\beta,\gamma)$ such that
$\nu$ is comparable with $\hat\eta(\alpha,\beta,\gamma,\delta)$.

Proof. Pick $t\in{\cal C}$ such that $\nu\in t$. Take  $\alpha$ and $\beta$ such that
$t$  is comparable with $\eta^*(\alpha,\beta)$. 

Case 1: $\nu\in\eta^*(\alpha,\beta)$.

Take $\gamma\leq k^{(\alpha,\beta)}$ and
$\delta\leq m^{(\alpha,\beta)}(\gamma)$ such that $\nu$ is comparable with $\eta^{(\alpha,\beta)}(\gamma,\delta)$.

Case 2: $\nu\notin\eta^*(\alpha,\beta)$.

If $\alpha = 0$ then let $\zeta= z(n+\beta+1)$ and if $\alpha>0$ then let
$\zeta= z(j^*(\alpha-1)+\beta+1)$. Let $\nu'=\nu\restr h(\zeta)$. Choose $\gamma\leq k^{(\alpha,\beta)}$ and
$\delta\leq m^{(\alpha,\beta)}(\gamma)$ such that $\nu'$ is comparable with
$\eta^{(\alpha,\beta)}(\gamma,\delta)$. Because the system $(j^{\alpha,\beta},
m^{\alpha,\beta)},k^{(\alpha,\beta)})$ is
$h$-tight above $\zeta$ we have $\eta^{(\alpha,\beta)}(\gamma,\delta)\leq \nu'$. Therefore
$\eta^{(\alpha,\beta)}(\gamma,\delta)\leq \nu$.

The Claim is established.

For each $n\in\omega$ define $\zeta(\alpha,\beta,\gamma)$ by the following recursive formulas:
$$\zeta(0,0,0)=0.$$
For $\alpha\leq k^*$ and $\beta\leq m^*(\alpha)$ and $\gamma<\hat k(\alpha,\beta)$ we have
$$\zeta(\alpha,\beta,\gamma+1)=\zeta(\alpha,\beta,\gamma)+1.$$
For $\alpha\leq k^*$ and $\beta< m^*(\alpha)$ we have
$$\zeta(\alpha,\beta+1,0)=\zeta(\alpha,\beta,\hat k(\alpha,\beta))+1.$$
For $\alpha< k^*$ we have
$$\zeta(\alpha+1,0,0)=\zeta(\alpha,m^*(\alpha),\hat k(\alpha,m^*(\alpha)))+1.$$

Define $\tilde\jmath(\zeta(\alpha,\beta,\gamma))=\hat\jmath(\alpha,\beta,\gamma)$, and
 $\tilde m(\zeta(\alpha,\beta,\gamma))=\hat m(\alpha,\beta,\gamma)$, and
 $\tilde k(\zeta(\alpha,\beta))=\hat k(\alpha,\beta)$, and
 $\tilde\eta(\zeta(\alpha,\beta,\gamma,\delta))=\hat\eta(\alpha,\beta,\gamma,\delta)$.

Claim 4. $(\tilde \jmath, \tilde k,\tilde m)$ is a $y$-bound system above $n$ and $H$ is a 
$(\tilde\jmath,\tilde k,\tilde m,\tilde\eta)$-squeezed tree.

Proof. By Claims 2 and 3.

The Lemma is established.

\proclaim Lemma 7.13.  Suppose $x\in{}^\omega(\omega-\{0\})$ is strictly increasing and suppose
that for each $n\in\omega$ we have that\/  $T_n$ is an $x$-squeezed tree.  Then there are
 $T^*$ and\/ $\langle\gamma_t\,\colon\allowbreak
t\in\omega\rangle$ an increasing sequence of integers such that\/
$T^*$ is an $x$-squeezed tree and $\gamma_0=0$ and\/
$(\forall t>0)\allowbreak(t<\gamma_t)$ and for every 
$f\in{}^{<\omega}\omega$ we have

\centerline{$(\forall t>0)(\exists s< t)(f\restr\gamma_t\in T_{\gamma_s})$ iff $f\in T^*$.}

\medskip

Proof: For each $n\in\omega$ choose $h_n\in{}^\omega\omega$ such that $T_n$ is $h_n$-tight.

We build as follows. Let $\gamma_0=0$.  Given $\gamma_t$, define $g_t(0)=\gamma_t$. For $0\leq s\leq t$ let
$g_t(s+1)=h_{\gamma_s}(g_t(s))$. Let $\gamma_{t+1}=g_t(t+1)$.

Let $T^*=\{\eta\in{}^{<\omega}\omega\,\colon\allowbreak (\forall t>0)(\exists s< t)
(\eta\restr \gamma_t\in T_{\gamma_s})\}$.

Now fix $n\in\omega$.  We build an $x$-bound system $(j,m,k)$ above $n$ and we build $\eta$ so that
$(j,m,k)$ and $\eta$ witness the fact that $T^*$ is $x$-squeezed.

For every $t\in\omega$ and $s\leq t$ choose an $h_{\gamma_s}$-tight $x$-bound system $(j^s_t,m^s_t,k^s_t)$ above
$g_t(s)$ along with $\eta^s_t$ such that
$T_{\gamma_s}$ is $(j^s_t,m^s_t,k^s_t,\eta^s_t)$-squeezed.

We define $\zeta$ such that for $\alpha\geq n$ and $\beta\leq \alpha$ and $\gamma\leq k^\alpha_\beta$ we have

\begin{itemize}

\item $\zeta(n,0,0)=0$, and

\item if $\gamma < k^\alpha_\beta$ then $\zeta(\alpha,\beta,\gamma+1)=\zeta(\alpha,\beta,\gamma)+1$, and

\item if $\beta< \alpha$ then $\zeta(\alpha,\beta+1,0)=\zeta(\alpha,\beta, k^\alpha_\beta)+1$, and

\item  if $\alpha\geq n$ then $\zeta(\alpha+1,0,0)=\zeta(\alpha,\alpha,k^\alpha_\alpha)+1$.

\end{itemize}

We define $(j,m,k)$ such that for every $\alpha\geq n$ and $\beta\leq \alpha$ and $\gamma\leq k^\alpha_\beta$ we have

\begin{itemize}

\item $j(\zeta(\alpha,\beta,\gamma))=j^\alpha_\beta(\gamma)$, and

\item $m(\zeta(\alpha,\beta,\gamma))=m^\alpha_\beta(\gamma)$, and

\item $k=k^\alpha_\beta$.

\end{itemize}

Claim 1. $(j,m,k)$ is an $x$-bound system above $n$.

Proof: Clause (1) of Definition 7.1 is immediate.

Clause (2) of Definition 7.1 holds because $j(0)=j^0_n(0)>x(g^0_n(0)+m^0_n(0)+1)\geq x(n+m(0)+1)$.  The first inequality holds because 
the system $(j^0_n,m^0_n,k^0_n)$ is above $g^0_n(0)$ and it satisfies clause (2) of Definition 7.1.

We have $j(\zeta(\alpha,\beta,\gamma+1))=j^\beta_\alpha(\gamma+1)>x(j^\beta_\alpha(\gamma)+m^\beta_\alpha(\gamma+1)+1)=
x(j(\zeta(\alpha,\beta,\gamma))+m(\zeta(\alpha,\beta,\gamma+1))+1)$.

We have $j(\zeta(\alpha,\beta+1,0))=j^{\beta+1}_\alpha(0)>x(g_\alpha(\beta+1)+m^{\beta+1}_\alpha(0)+1)\geq
x(h_{\gamma_\beta}(g_\alpha(\beta))+m^{\beta+1}_\alpha(0)+1)\geq
x(j^\beta_\alpha(k^\beta_\alpha)+m^{\beta+1}_\alpha(0)+1)=x(j(\zeta(\alpha,\beta,k^\beta_\alpha))+m(\zeta(\alpha,\beta+1,0))+1)$.

The first inequality is clause (2) of Definition 7.1 applied to the system $(j^{\beta+1}_\alpha,\allowbreak
m^{\beta+1}_\alpha,k^{\beta+1}_\alpha)$.  The second inequality is by the definition
 of $g_\alpha$. The third inequality is because the system
$(j^\beta_\alpha,m^\beta_\alpha,k^\beta_\alpha)$ is $h_{\gamma_\beta}$-tight above $g_\alpha(\beta)$.

The Claim is established.

We define $\eta$ such that for every $\alpha\geq n$ and $\beta\leq \alpha$ and $\gamma\leq k^\alpha_\beta$ and $\delta\leq m^\beta_\alpha(\gamma)$
we have $\eta(\zeta(\alpha,\beta,\gamma),\delta)=\eta^\beta_\alpha(\gamma,\delta)$.

Claim 3: $T^*$  is a $(j,k,m,\eta)$-squeezed tree. 

Proof: It is straightforward to verify that $T^*$ is a tree and that clause (1) and clause (2) of Definition 7.2 hold.

To verify clause (3), suppose we have $\nu\in T^*$. We show that $\nu$ is comparable to some $\eta(l,i)$
with $(l,i)\in{\rm dom}(\eta)$. Choose $\nu'\in T^*$ such that $\nu\leq \nu'$ and ${\rm lh}(\nu')\geq
\gamma_{n+1}$. It suffices to show that $\nu'$ is comparable with some $\eta(l,i)$ with $(l,i)\in {\rm dom}(\eta)$.
Because $\nu'\in T^*$ we may choose $s\leq n$ such that $\nu'\restr\gamma_{n+1}\in T_{\gamma_s}$.
We may select $(l,i)\in{\rm dom}(\eta^s_n)$ such that $\nu'\restr\gamma_{n+1}$ is comparable with $\eta^s_n(l,i)$.
We have $\eta^n_s(l,i)=\eta(\zeta(n,s,l),i)$, so ${\rm lh}(\eta^s_n(l,i))=j(\zeta(n,s,l))=j^s_n(l)\leq
h_{\gamma_s}(g_n(s))=g_n(s+1)\leq\gamma_{n+1}$. Therefore $\eta(\zeta(n,s,l),i)\leq\nu'\restr\gamma_{n+1}$ and therefore
$\eta(\zeta(n,s,l),i)$ is comparable with $\nu'$.

The Claim and the Lemma are established.

The following Lemma is the analogue of Lemma 4.7.  The fact that the Lemma is stronger reflects the fact that
Lemma 7.10 holds.

\proclaim Lemma 7.14.  Suppose $x\in{}^\omega(\omega-\{0\})$ is strictly increasing, and suppose
 $P$ is a forcing notion such that\/ {\rm $V[G_P]\models``$for all countable $X\subseteq V$ there is
a countable $Y\in V$ such that $X\subseteq Y$ and $\langle T_n\,\colon\allowbreak
n\in\omega\rangle$ is a sequence of  $x$-squeezed trees and $(\forall n\in\omega)\allowbreak(T_n\in V)$.''}
Then\/ {\rm $V[G_P]\models``$there is a strictly increasing sequence of integers $\langle m_i\,\colon\allowbreak
i\in\omega\rangle$ and an $x$-squeezed tree $T^*\in V$ such that $m_0=0$ and
$(\forall i>0)\allowbreak(m_i>i)$ and
for every $\eta\in{}^{<\omega}\omega$, if
$(\forall i>0)\allowbreak(\exists j<i)\allowbreak(\eta\restr m_{i+1}\in T_{m_j})$ then $\eta\in T^*$.''}

Proof: Work in $V[G_P]$. Let $b\in V$ be a countable set such that
$\{T_n\,\colon\allowbreak n\in\omega\}\subseteq b\in V$ and $(\forall x\in b)\allowbreak(x$ is an $x$-squeezed tree).
Let $\langle S_n\,\colon\allowbreak n\in\omega\rangle\in V$ enumerate $b$ with infinitely many repetitions such that
$S_0=T_0$.
Build $\langle S'_n\,\colon\allowbreak n\in\omega\rangle$ by setting $S'_0=S_0$ and for every $n>0$ set
$S'_n=S_n\cup S'_{n-1}$. Build $h$ mapping $\omega$ into $\omega$  by setting $h(0)=0$ and for every
$n>0$ set $h(n)$ equal to the least integer  $m>n$ and $T_n=S_m$.

Using Lemma 7.13, take $ T^*\in V$  an $x$-squeezed tree and $\langle k_i\,\colon\allowbreak i\in\omega\rangle$ such that
for every $\eta\in{}^{<\omega}\omega$ we have $\eta\in T^*$ iff
$(\forall n>0)\allowbreak(\exists i<n)\allowbreak(\eta\restr k_n\in S'_{k_i})$.

Build $\langle n'_i\,\colon\allowbreak i\in\omega\rangle$ an increasing sequence of integers such that
$n'_0=0$ and $n'_1>k_1$ and for every $i\in\omega$ we have $h(n'_i)<n'_{i+1}$ and 

(*)  $(\exists t\in\omega)\allowbreak
(n'_i<k_t<n'_{i+1})$.

For every $i\in\omega$ let $m_i=n'_{2i}$.

Fix $\eta\in{}^{<\omega}\omega$ such that $(\forall i>0)(\exists j<i)(\eta\restr m_{i+1}\in T_{m_j})$. To establish the Lemma, it suffices to show $\eta\in T^*$.
By choice of $T^*$, it suffices to show $(\forall n>0)\allowbreak(\exists i<n)\allowbreak(\eta\restr k_{n}\in S'_{k_i})$.

Claim 1. $(\forall i>0)(\exists j<i)(\eta\restr n'_{i+1}\in S'_{n'_j})$.

Proof: The proof breaks into two cases.

Case 1: $i<4$.

We have $n'_{i+1}\leq n'_4= m_2$, and $\eta\restr m_2\in T_0$, so
$\eta\restr n'_{i+1}\in S_0'$.

Case 2: $i\geq 4$.

Fix $i^*>0$ such that $2i^*\leq i \leq 2i^*+1$.

We may fix $j^*<i^*$ such that $\eta\restr m_{i^*+1}\in T_{m_{j^*}}$.  

Now, we have 

(*)  $i+1\leq 2i^*+2$ so

(**) $ n'_{i+1}\leq  m_{i^*+1}$. 

 We also have 

(***) $\eta\restr m_{i^*+1}\in T_{m_{j^*}}\subseteq S'_{h(m_{j^*})}$.

By (**) and (***) we have

(****)  $\eta\restr n'_{i+1}\in S'_{h(m_{j^*})}$.

Note that

(*****) $h(m_{j^*})=h(n'_{2j^*})\leq n'_{2j^*+1}\leq n'_{2i^*-1}\leq n'_{i-1}$.

By (****) and (*****) we have $\eta\restr n'_{i+1}\in S'_{h(m_{j^*})}\subseteq S'_{n'_{i-1}}$.

The Claim is established.

To complete the proof of the Lemma, suppose $i>0$. We must show that there is $t<i$ such that
$\eta\restr k_{i}\in S'_{k_t}$.

Case 1: $k_{i-1}<n'_0$.

By (*) we have $n'_1\geq k_{i}$. By Claim 1 we have $\eta\restr n'_1\in S_0$.  Hence
$\eta\restr k_{i}\in S_0$.

Case 2:  $n'_0\leq k_{i-1}$.

By (*) we know that there is at most one element of $\{n'_j\,\colon\allowbreak j\in\omega\}$ strictly between
$k_{i-1}$ and $k_{i}$. Hence we may
fix $j>0$ such that $n'_{j-1}\leq k_{i-1}<k_{i}\leq n'_{j+1}$.  If $\eta\restr n'_{j+1}\in S_0$ then
$\eta\restr k_{i}\in S_0$ and we are done, so assume otherwise.  By Claim 1 we may fix $m<j$ such that $\eta\restr n'_{j+1}\in S'_{n'_m}$.
We have $\eta\restr k_{i}\in S'_{n'_m}\subseteq S'_{n'_{j-1}}\subseteq S'_{k_{i-1}}$ and again we are done.

The Lemma is established.

The following Lemma is the analogue of Lemma 4.9.

\proclaim Lemma 7.15.  Suppose
$y\in{}^\omega(\omega-\{0\})$ is strictly increasing and\/
 $P$ is a forcing notion such that\/ {\rm $V[G_P]\models``$for all countable $X\subseteq V$ there is
a countable $Y\in V$ such that $X\subseteq Y$ and
$\langle T_n\,\colon\allowbreak
n\in\omega\rangle$ is a sequence of  $y$-squeezed trees, each of which is in $V$.''}
 Then in $V[G_P] $
there is a $y$-squeezed tree
$T^*\in V$ such that   for every
$n\in\omega$ and every $j\in\omega$ and every $g\in[T_j]$
there is $k\in\omega$ such that for every $\eta\in T_j$ extending $g\restr k$, if 
 $\eta\restr k\in  T^*$ then
$\eta\in  T^*$.

Proof:  In $V[G]$, build a sequence of $y$-squeezed trees  $\langle T'_j\,\colon j\in\omega\rangle$, each in $V$, such that $T'_0=T_0$ and for every
$j\in\omega$ we have $T'_{j+1}=T'_j\cup T_{j+1}$. 
By Lemma 7.14 we may find an increasing sequence of integers 
$\langle k_n\,\colon\allowbreak n\in\omega\rangle$ and a $y$-squeezed tree 
$T^*\in V$ such that 
$(\forall n>0)\allowbreak(k_n>n)$ and
for every $\eta\in{}^{<\omega}\omega$ we have

\medskip

{\centerline{$(\forall n>0)(\exists i< n)(
\eta\restr k_n\in T'_{k_i})$ iff $\eta\in T^*$.}}

\medskip

 Fix  $j\in\omega$ and $g\in [T_j]$.
Let $k={\rm max}\{k_{j'}\colon\allowbreak j'\leq j\}$. 
 Fix $\eta\in T_j$ extending $g\restr k$ and assume
$\eta\restr k\in T^*$.  It suffices to show that $\eta\in T^*$. 
If $j=0$ then $\eta\in T_0=T'_0\subseteq T^*$.
Therefore, we assume that $j>0$.  It suffices to show that

\medskip

\centerline{ $(\forall i>0)\allowbreak(\exists i'< i)
\allowbreak(\eta\restr k_i\in T'_{k_{i'}})$. }

\medskip

 Towards this end, fix $i>0$.

Case 1: $i\leq j$.

Because $k_i\leq k$ we have $\eta\restr k_i\in T^*$. Therefore we may take $i'< i$ such that
$\eta\restr k_{i}\in T'_{k_{i'}}$.  

Case 2: $0<j<i$.

Because $\eta\in T_j$ we  have $\eta\restr i\in T_j\subseteq T'_{j+1}\subseteq T'_{k_j}$.

The Lemma is established.

The following is the analogue of Lemma 4.8.

\proclaim Lemma 7.16.  Suppose   $y\in{}^\omega(\omega-\{0\})$ is strictly increasing
and\/ $\zeta\in{}^\omega\omega$  and $\langle T_n\,\colon\allowbreak
n\in\omega\rangle$ is a sequence of  $y$-squeezed trees.
 Then 
there is a $y$-squeezed tree $T^*$ and a sequence of integers
$\langle m_i\,\colon\allowbreak i\in\omega\rangle$
such that $\zeta\in[T^*]$ and $T\subseteq T^*$ and
for every  $i\in\omega$ and every $j>m_i$ and every $\nu\in T_{m_i}$ extending $\zeta\restr j$
we have  $\nu\in T^*$.

Proof.
 Define $\langle T'_k\,\colon\allowbreak k\in\omega\rangle $ by setting $T'_0=
 T_0\cup\{\zeta\restr n\,\colon\allowbreak n\in\omega\}$ and for every $k\in\omega$
set $T'_{k+1}= T'_k\cup T_{k+1}$.

By Lemma 7.14 we may choose $T^*$ a $y$-squeezed tree and
$\langle m_i\,\colon\allowbreak i\in\omega\rangle$ an increasing sequence of integers
such that

\medskip

{\centerline{$(\forall g\in{}^\omega\omega)\allowbreak((\forall n>0)\allowbreak
(\exists i< n)\allowbreak(g\restr m_n\in T'_{m_i})$ implies $g\in[T^*])$.}}

\medskip

Now suppose that $\eta\in T$ and $i\in\omega$ and
${\rm length}(\eta)\geq m_i$ and
 $\nu$ extends $\eta$
and $\nu\in T_{m_i}$.  We show $\nu\in T^*$.

Choose $h\in[T_{m_i}]$ extending $\nu$.  It suffices to show that $h\in[T^*]$.
Therefore it suffices to show that $(\forall k>0)\allowbreak(\exists
j< k)\allowbreak(h\restr m_k\in T'_{m_j})$.

Fix $k\in\omega$. 
If $i< k$ then because $h\in[T_{m_i}]$ we have that
$h\restr m_k\in T_{m_i}\subseteq T'_{m_i}$ and we are done.  If instead $k\leq i$ then
$h\restr m_k=\eta\restr m_k\in T'_{0}$ and again we are done.

The Lemma is established.

The following Theorem is [13, Theorem VI.1.12] for the case of the $P$-point property.  Rather than simply referring to the proof of Theorem 4.10, we give the 
complete argument to demonstrate the simplifications afforded us by the fact that Lemma 7.10 holds.

\proclaim Theorem 7.17.  Suppose $\langle P_\eta\,\colon\allowbreak\eta\leq\kappa\rangle$ is a countable support iteration based on\/
 $\langle Q_\eta\,\colon\allowbreak\eta<\kappa\rangle$ and suppose\/ {\rm $(\forall\eta<\kappa)\allowbreak({\bf 1}\forces_{P_\eta}``Q_\eta$ 
is proper and has the $P$-point property'').} 
Then $P_\kappa$ has the $P$-point property.

Proof: By induction on $\kappa$.  No counterexample can first appear 
at a stage of uncountable cofinality, and the successor case is
easily handled using Lemma 7.12, so we may assume $\kappa =\omega$.

Suppose $\lambda$ is a sufficiently large regular cardinal  
and   $\zeta$ is a $P_\omega$-name 
and $y\in
{}^\omega(\omega-\{0\})$ is strictly increasing and\/ {\rm
 ${\bf 1}\forces_{P_\omega}``\zeta\in{}^\omega\omega$.''}
Suppose $N$ is a countable elementary submodel of $H_\lambda$ and $\{P_\omega,  y, \zeta
\}\in N$.

Fix $p$ in $P_\omega\cap N$.

Fix $p'\in N$ and
 $\langle (p_n,\zeta_n)\,\colon\allowbreak n\in\omega\rangle\in N$ 
as in Lemma 3.1.

Let $\langle T'_{j}\,\colon\allowbreak j\in\omega\rangle$ list all  $T'\in N$ such that 
 we have that $T'$ is a
$y$-squeezed tree, with infinitely many repetitions.

Working in $V[G_{P_\omega}]$, use Lemma 7.15 to choose $T^*\in V$ 
a $y$-squeezed tree such that 
for every $n\in\omega$ and every $j\in\omega$ and every $g\in[T'_j]$ there exists
$k\in\omega$ such that for every $\eta\in T'_j$ extending $g\restr k$, if
$\eta\restr k\in T^*$ then $\eta\in T^*$.

In the preceding paragraph, we worked in $V[G_{P_\omega}]$ so that the brackets about $T'_j$ would be interpreted in $V[G_{P_\omega}]$; i.e.,
$g$ need not be in $V$.

Claim 1.  We may be build $\langle r_n\,\colon\allowbreak n\in\omega\rangle$ such that
 for every $n\in\omega$ we have that the following hold:

(1) $r_n\in P_{n}$ is $N$-generic, and

(2) $r_{n+1}\restr n=r_n$, and

(3) $r_n\forces``\zeta_n\in[T^*]$,'' and

(4) $r_n\leq p'\restr n$.''

Proof: By induction on $n$.  For $n=0$ we have nothing to prove. Suppose we have $r_n$.

Let $F_0$ and $F_2$  be $P_{n}$-names  such that

\medskip

(*) {
 ${\bf 1}\forces``F_0$ and
$F_2$ are functions and each of whose domains
 is equal

\centerline{to $Q_n$, such that}
$(\forall q'\in Q_n )\allowbreak (F_0(q')$ is a $y$-squeezed tree}

\centerline{and
$F_2(q')\leq q'$
and $F_2(q')\forces`\zeta_{n+1}\in[ F_0(q')]$').''}

\medskip

 We may assume that the names $F_0$ and
$F_2$  are in $N$. Notice that $F_0$ and $F_2$ depend on $n$,
although this dependence is suppressed in our notation.

Working in $V[G_{P_{n}}]$, use Lemma 7.16 to choose 
$ T^\#_n$ a $y$-squeezed tree in $V$ and
$\langle k_i\,\colon\allowbreak i\in\omega\rangle$ an increasing sequence of
integers (this sequence depends on $n$ but this fact is suppressed in our notation)
such that 
$\zeta_n\in[ T^\#_n]$
and for every $\eta$  and every $i\in\omega$ and every
$\nu\in F_0(p_{k_i}(n))$,
 if 
$\eta$ is a proper initial segment of $\zeta_n$ and ${\rm length}(\eta)\geq  k_i$ and 
$\nu$ extends $\eta$, then $\nu\in T^\#_n$. 

\medskip

We may assume the $P_{n}$-name $ T^\#_n$ is in $N$.

 Using the induction hypothesis and Lemma 7.12,  fix $\tilde T_n\in V$ a $y$-squeezed tree such that
$T^\#_n\subseteq \tilde T_n$.

Because  $\tilde T_n$ is a $P_{n}$-name in $N$ forced to be in $V$,
 we conclude that by the $N$-genericity of $r_n$ that

\medskip

{\centerline{$r_n\forces``\tilde T_n\in N$.''}}

\medskip

Therefore there is a $P_{n}$-name $m$ such that

\medskip

{\centerline{$r_n\forces``\tilde T_n=T_m'$ and $m>n$.''}}

\medskip

Because $T^*$ was chosen as in the conclusion of Lemma 7.15, we may choose
 $k$ to be a $P_{n}$-name for an integer such that

\medskip

(**) $r_n\forces``(\forall\eta\in \tilde T_n)($if $\eta$ extends 
$\zeta_n\restr k$ and $\eta\restr k\in  T^*$ 
then $\eta\in  T^*)$.''

\medskip

Choose $j$  to be a $P_{n}$-name such that
$r_n\forces``k_j\geq k$.''

Subclaim 1. 
$r_n\forces``F_2(p_{k_j}(n))\forces
`\zeta_{n+1}\in[\tilde T_n]$.'\thinspace''

Proof.  It suffices to show

\medskip 

$r_n\forces``F_2(p_{k_j}(n))\forces
`(\forall j'>j)(\zeta\restr k_{j'}\in\tilde T_n)$.'\thinspace''

\medskip

Fix $j'$ a $P_{{n+1}}$-name for an integer such that

\medskip

$r_n\forces``F_2(p_{k_j}(n))\forces
`j'>j$,'\thinspace''

\medskip

By the definition of $\langle p_i\,\colon\allowbreak i\in\omega\rangle$ we have

\medskip

(***) $r_n\forces``p_{k_j}(n)\forces`\zeta_n\restr k_j=\zeta_{n+1}\restr k_j$.'\thinspace''

\medskip

By (*) we have

\medskip

(****) $r_n\forces``F_2(p_{k_j}(n))\forces`
\zeta_{n+1}\in [F_0(p_{k_j}(n))]$.'\thinspace''

\medskip

Combining (***), (****), and the definition of $\tilde T_n$, we have that

\medskip

{\centerline{$r_n\forces``F_2(p_{k_j}(n))\forces
`\zeta_{n+1}\restr k_{j'}\in \tilde T_n$.'\thinspace''}}

The Subclaim is established.

Subclaim 2. $r_n\forces``F_2(p_{k_j}(n))\forces`\zeta_{n+1}\in[ T^*]$.'\thinspace''

Proof: By (**) we have

\medskip

($\dag$) $r_n\forces``(\forall\eta\in \tilde T_n)\allowbreak(\eta\restr k_j\in  T^*$ implies
$\eta\in T^*)$.''

\medskip

Work in $V[G_{P_{n}}]$ with $r_n\in G_{P_{n}}$.
Fix $\eta\in\tilde T_n$ and suppose $F_2(p_{k_j}(n))\forces``
\eta$ is an initial segment of $\zeta_{n+1}$ and ${\rm lh}(\eta)\geq k_j$.''
To establish the Subclaim it suffices to show

($\#)$ $F_2(p_{k_j}(n))\forces``\eta\in T^*$.''

By the definition of $\langle p_i\,\colon\allowbreak i\in\omega\rangle$ we have

\medskip

{\centerline{$p_{k_j}(n)\forces``\eta\restr k_j=\zeta_n\restr k_j$.''}}

\medskip

Hence by the fact that Claim 1 holds for the integer $n$ we have

\medskip

($\dag\dag$) $p_{k_j}(n)\forces``\eta\restr k_j\in  T^*$.''

\medskip

By Subclaim 1, ($\dag$), ($\dag\dag$), and the fact that
$F_2(p_{k_j}(n))\leq
p_{k_j}(n)$ we obtain

\medskip

{\centerline{$F_2(p_{k_j}(n))\forces``\eta\in  T^*$.''}}

\medskip

Subclaim 2 is established.

To complete the induction establishing Claim 1, we  take
$r_{n+1}\in P_{{n+1}}$ such that $r_{n+1}\restr n=r_n$ and
$r_{n+1}$ is $N$-generic and $r_n\forces``r_{n+1}(n)\leq
F_2(p_{k_j}(n))$.''

Claim 1 is established.

Define $q$  by

\medskip

\centerline{$q=\bigcup\{r_n\,\colon\allowbreak n\in\omega\}$.}

\medskip

We have $q\leq p$ and $q\forces``\zeta\in[T^*]$.''

The Theorem is established.

\section{On adding no Cohen reals}

In [13, Conclusion VI.2.13D(1)], Shelah states that a countable support iteration of proper forcings, each of which adds no Cohen
reals, either adds no Cohen reals or adds a dominating real.  However, according to Jakob Kellner,
Shelah has stated that this is an error, and the result holds only at limit stages.
In this section, we prove the limit case.

\proclaim Definition 8.1. A nowhere dense tree $T\subseteq{}^{<\omega}\omega$ is a non-empty tree such that
for every $\eta\in T$ there is some $\nu$ extending $\eta$ such that $\nu\notin T$. A perfect tree $T\subseteq{}^{<\omega}\omega$ is a 
non-empty tree such that
  for every $\eta\in T$, the set of successors of $\eta$ in $T$ is not linearly ordered.

\proclaim Lemma 8.2.  $P$ does not add any Cohen reals iff\/ {\rm ${\bf 1}\forces_P``(\forall f\in{}^\omega\omega)\allowbreak(\exists H\in V)\allowbreak(H$ is a nowhere dense perfect tree and $f\in[H])$.''}

Proof: This is a tautological consequence of the definition of Cohen real.

\proclaim Lemma 8.3.  Suppose ${\rm cf}(\kappa)=\omega$ and\/ $\langle\alpha_n\,\colon n\in\omega\rangle$ is an increasing sequence of ordinals cofinal in
$\kappa$ such that $\alpha_0=0$.  Suppose
$\langle P_\eta\,\colon\allowbreak \eta\leq\kappa\rangle$ is a countable support forcing iteration based
on $\langle Q_\eta\,\colon\allowbreak\eta<\kappa\rangle$ and
for every $\eta<\kappa$ we have {\rm $V[G_{P_\eta}]\models``Q_\eta$ is proper.''}
Suppose $p\in P_\kappa$ and\/ {\rm $p\forces``f\in{}^\omega\omega$.''}
Then there are $p'\leq p$ and $\langle \eta_n\,\colon\allowbreak n\in\omega\rangle$ such that
for every $n\in\omega$ we have that $\eta_n$ is a $P_{\alpha_n}$-name and\/ {\rm
$p'\restr\alpha_n\forces``p'\restr[\alpha_n,\kappa)\forces`\eta_n= f\restr n$.'\thinspace''}

Proof: Let $\lambda$ be a sufficiently large regular cardinal and $N$ a countable elementary
substructure of $H_\lambda$ such that $\{P_\kappa,p,f,\langle \alpha_n\,\colon\allowbreak
n\in\omega\rangle\}\in N$.

Using the Proper Iteration Lemma, build $\langle(p_n,q_n,\eta_n)\,\colon\allowbreak
n\in\omega\rangle$ such that $q_0=p$ and for every $n\in\omega$ we have the following:

(1) $p_n\in P_{\alpha_n}$ is $N$-generic and $\eta_n$ is a $P_{\alpha_n}$-name, and

(2)  $p_n\forces``q_{n+1}\leq q_n\restr[\alpha_n,\kappa)$ and $q_{n+1}\in N[G_{P_{\alpha_n}}]$ and $q_{n+1}\forces`\eta_n=f\restr n$,'\thinspace'' and

(3) $p_{n+1}\restr \alpha_n = p_n$ and $p_n\forces``p_{n+1}\restr[\alpha_n,\alpha_{n+1})\leq q_{n+1}\restr\alpha_{n+1}$.''

Letting $p'=\bigcup\{p_n\,\colon n\in\omega\}$ establishes the Lemma.

The proof of the following Theorem is [13, proofs of Claims VI.2.5(2) and VI.2.13C]. 

\proclaim Theorem 8.4.  Suppose ${\rm cf}(\kappa)=\omega$ and
$\langle P_\eta\,\colon\allowbreak\eta\leq\kappa\rangle$ is a countable
support forcing iteration based on $\langle Q_\eta\,\colon\allowbreak
\eta<\kappa\rangle$ and for every $\eta<\kappa$ we have {\rm
$V[G_{P_\eta}]\models``Q_\eta$ is proper''} and $P_\eta$  does not add any Cohen reals.
  Suppose $P_\kappa $ does not add any dominating reals. Then $P_\kappa$ does not add any Cohen reals.

Proof. Fix $\langle\alpha_n\,\colon\allowbreak n\in\omega\rangle$ cofinal in $\kappa$ with $\alpha_0=0$. Also in $V[G_{P_\kappa}]$ fix
$f\in{}^\omega\omega$. Suppose $p\in P_\kappa$ and $p\forces``f$ is a Cohen real.''

Let $p'\leq p$ and $\langle\eta_n\,\colon n\in\omega\rangle $ be as in Lemma 8.3.

Working in $V[G_{P_\kappa}]$ with $p'\in G_{P_\kappa}$, 
let $\langle T_n\,\colon\allowbreak n\in\omega\rangle$ be a sequence of nowhere dense perfect trees such that
$(\forall n\in\omega)\allowbreak(T_n\in V$ and $\eta_n\in  T_n)$.

 Let $B\in V$ be a countable set of nowhere dense perfect trees  such that for every $n\in\omega$ we have
$T_n\in B$.
Let $\langle S_n\,\colon\allowbreak n\in\omega\rangle\in V$ enumerate  $B$ with infinitely many repetitions such that $T_0=S_0$.

Build inductively $\langle S'_n\,\colon\allowbreak n\in\omega\rangle$ such that
$S'_{n+1}=S_{n+1}\cup S'_n$ and $S_0'=S_0$.

Define $h\in{}^\omega\omega$ by setting $h(k) $ equal to the least $m>k$ such that $T_{k}\subseteq S'_m$, for every $k\in \omega$.
Because $P_\kappa$ adds no dominating reals we may choose $g\in{}^\omega\omega\cap V$  and $A\subseteq \omega$ such that 
$A=\{n\in \omega\,\colon\allowbreak g(n)>h(n)\}$ and $A$ is infinite.

Choose $\langle k_i\,\colon\allowbreak i\in\omega\rangle\in V$ an increasing sequence of integers
as follows.
Let $k_0=0$. Given $k_n$, choose $k_{n+1}\geq{\rm max}(k_n+1,2)$ such that $(\forall\nu\in{}^{\leq k_n}k_n)\allowbreak
[(\exists\nu'\in{}^{k_{n+1}}\omega$ extending $\nu)\allowbreak(\forall i\leq k_n)\allowbreak(\nu'\notin S'_{g(i)})$ and
$(\forall i\leq k_n)\allowbreak(\exists \nu_1\in S'_{g(i)})\allowbreak(\exists \nu_2\in S'_{g(i)})\allowbreak(\nu_1$ and $\nu_2$ 
are distinct extensions of $\nu$ and
${\rm lh}(\nu_1)={\rm lh}(\nu_2)=k_{n+1})]$.

Let $T^0=\{\eta\in{}^{<\omega}\omega\,\colon\allowbreak(\exists s\in\omega)\allowbreak(\exists j\in\omega)\allowbreak
(k_{2s}\leq j<k_{2s+1}$ and $\eta\restr j\in S'_0$ and $\eta\in S'_{g(j)})\}$.

Let $T^1=\{\eta\in{}^{<\omega}\omega\,\colon\allowbreak(\exists s\in\omega)\allowbreak(\exists j\in\omega)\allowbreak
(k_{2s+1}\leq j<k_{2s+2}$ and $\eta\restr j\in S'_0$ and $\eta\in S'_{g(j)})\}$.

Claim 1: $T^0$ is a nowhere dense  tree.

Proof. Suppose $\eta\in T^0$. Choose $s$ and $j$ witnessing this. Also take $n\geq s$ so large that
$\eta\in{}^{\leq k_{2n}}k_{2n}$.

We  choose $\nu$ extending $\eta$ such that ${\rm lh}(\nu)=k_{2n+2}$ and
$(\forall i\leq k_{2n+1})\allowbreak(\nu\notin S'_{g(i)})$. In particular we have
$\nu\notin S'_0$.  We show that $\nu\notin T^0$.
So suppose, towards a contradiction, that $s'\in\omega$ and
$j'\in\omega$ and $k_{2s'}\leq j'<k_{2s'+1}$ and $\nu\restr j'\in S'_0$ and $\nu\in S'_{g(j')}$. Because $\nu\in S'_{g(j')}$ we know $j'\geq k_{2n+1}$.
Necessarily, then, $j'\geq k_{2n+2}$. But then $\nu=\nu\restr j'\in S'_0$.  
This  contradiction establishes the 
Claim.

Claim 2. $T^0$ is a perfect tree.

Proof: Given $\eta\in T^0$, let $s\in\omega$ and $j\in\omega$ be witnesses. 

Case 1: ${\rm lh}(\eta)\geq j$.

Let $\nu$ and $\nu'$ be incomparable elements of $S'_{g(j)}$
extending $\eta$.  We have that $\nu$ and $\nu'$ are in $T^0$; 
this is witnessed by the integers $s$ and $j$.

Case 2: ${\rm lh}(\eta)<j$.

Take $\nu$ and $\nu'$ distinct extensions of $\eta$ such that $\nu\in S_0$ and $\nu'\in S_0$ and
${\rm lh}(\nu)={\rm lh}(\nu')=j$. We have $\nu\in S'_{g(j)}$ and $\nu'\in S'_{g(j)}$ because
$S_0\subseteq S'_{g(j)}$.  We have that $\nu$ and $\nu'$ are in $T^0$;
 this is witnessed by the integers $s$ and $j$.

Claim 3: $T^1$ is a nowhere dense perfect tree.

Proof: Similar to Claims 1 and 2.

Let $B_0=\bigcup\{[k_{2i},k_{2i+1})\,\colon\allowbreak i\in\omega\}$ and let $B_1=\bigcup\{[k_{2i+1},k_{2i+2})\,\colon\allowbreak i\in\omega\}$.

Claim 4:  $(\forall n\in A\cap B_0)\allowbreak(\eta_{n}\in T^0)$.
  
Proof: Given  $n\in A\cap B_0$ choose $s\in\omega$ such that
$k_{2s}\leq n< k_{2s+1}$. We have
$\eta_{n}\in T_{n}\subseteq S'_{h(n)}\subseteq S'_{g(n)}$ and $\eta_n\in S'_0$. Hence
$\eta_{n}\in T^0$.  The Claim is established.

Claim 5:  $(\forall n\in A\cap B_1)\allowbreak(\eta_{n}\in T^1)$.

Proof: Similar to Claim 4.

We have that $T^0$ and $T^1$ are elements of $V$. Furthermore, if $A\cap B_0$ is infinite, we have by Claim 4  that 
for infinitely many $n$ we have $\eta_{n}\in T^0$ and hence $f\in[T^0]$.
Otherwise by Claim 5 it follows  that for infinitely many $n$ we have $\eta_n\in T^1$ 
and hence $f\in[T^1]$.  The Theorem is established.

\section{On not adding reals not belonging to any closed null sets of $V$}

In this section we give Shelah's proof that the property $``P$ does not add any real not belonging to any closed set of measure zero of the ground model'' is preserved at limit stages by countable support iterations of proper forcings assuming the iteration does not add dominating reals.

\proclaim Theorem 9.1.    Suppose $\langle P_\eta\,\colon\allowbreak\eta\leq\kappa\rangle$ is a countable support iteration based on $\langle Q_\eta\,\colon\allowbreak\eta<\kappa\rangle$ and suppose\/ $\kappa$ is a limit ordinal and\/ {\rm $(\forall\eta<\kappa)\allowbreak(P_\eta$ 
 does not add reals not in any closed measure zero set of $V$).}   Suppose also that $P_\kappa$ does not add any dominating reals.
Then $P_\kappa$ does not add any real not in any closed measure zero set of $V$.

Proof:  Repeat the proof of Theorem 8.3 with ``nowhere dense perfect tree'' replace by ``perfect tree with Lebesgue measure zero''
throughout, and choosing $k_{n+1}$ so large that $2^{-k_{n+1}}\cdot\vert\{\nu\in {}^{k_{n+1}}\omega\,\colon\allowbreak
(\exists i\leq k_n)\allowbreak(\nu\in S'_{g(i)})\}\vert < 1/n$.

\vfill\eject

\noindent{\bf References}

\bigskip

[1] Abraham, U., Proper Forcing, in {\bf Handbook of Set Theory, vol.~1}, Springer, 2010

\medskip 

[2] Eisworth, T., CH and first countable, countably compact spaces, Topology Appl. 109, no. 1, 55-73 (2001)

\medskip

[3]  Goldstern, M., Tools for Your Forcing Construction, Set theory of the reals (Haim Judah, editor), Israel Mathematical Conference Proceedings, vol. 6, American Mathematical Society, pp. 305--360. (1993)

\medskip

[4]  Goldstern, M. and J. Kellner, New reals: Can live with them, can live without them, Math. Log. Quart. 52, No. 2, pp.~115--124 (2006)

 \medskip

[5] Kellner J., and S. Shelah, Preserving preservation, Journal of Symbolic Logic, vol 70, pp.~914--945 (2005)

\medskip

[6] Schlindwein, C., Consistency of Suslin's hypothesis, a non-special Aronszajn tree, and GCH, Journal of Symbolic Logic, vol.~59, pp.~1--29, 1994

\medskip

[7] Schlindwein, C., Suslin's hypothesis does not imply stationary antichains, Annals of Pure and Applied Logic, vol.~64, pp.~153--167, 1993

\medskip

[8] Schlindwein, C., Special non-special $\aleph_1$-trees, Set theory and its Applications,  J.~Steprans and S.~Watson (eds.), Lecture
Notes
in Mathematics, vol.~1401, Springer-Verlag, 1989

\medskip

[9] Schlindwein, C., Shelah's work on non-semi-proper iterations, II, Journal of Symbolic Logic,
vol.~66 (4), pp.~1865 -- 1883, 2001

\medskip

[10] Schlindwein, C., SH plus CH does not imply stationary antichains, Annals of Pure and Applied Logic, vol.~124, pp.~233--265, 2003

\medskip

[11] Schlindwein, C., Shelah's work on non-semi-proper iterations, I, Archive for Mathematical Logic, vol.~47 (6), pp.~579 -- 606, 2008

\medskip

[12] Schlindwein, C., Understanding preservation theorems, II,  Mathematical Logic Quarterly,
vol.~56, pp.~549--560, 2010

\medskip

[13] Shelah, S., {\bf Proper and Improper Forcing}, Perspectives in Mathematical Logic, Springer, Berlin, 1998

\end{document}